\documentclass[10pt]{article}

\usepackage[T1]{fontenc}

\usepackage[english]{babel}

\usepackage{amsmath,amsthm,amssymb}
\usepackage{graphicx}
\usepackage[all]{xy}
\usepackage{psfrag}

\usepackage{color}
\usepackage{epsfig}
\usepackage{latexsym}

\newtheorem{thm}{Theorem}[section]
\newtheorem{cor}[thm]{Corollary}
\newtheorem{lem}[thm]{Lemma}
\newtheorem{prop}[thm]{Proposition}

\newtheorem{rem}[thm]{Remark}
\newtheorem{rems}[thm]{Remarks}
\newtheorem*{acknowledgement}{Acknowledgements}

\numberwithin{figure}{section}

\newcommand{\pnc}{\mathbb{C} \mathbb{P}^n}
\newcommand{\cc}{\mathbb{C}}
\newcommand{\rr}{\mathbb{R}}
\newcommand{\zz}{\mathbb{Z}}

\newcommand{\TT}{\mathbb{T}}
\newcommand{\grad}{\mbox{grad}}
\newcommand{\cotg}{T^{\ast}V}
\newcommand{\ctg}{T^{\ast}}
\newcommand{\tH}{\widetilde{H}}
\newcommand{\tPhi}{\widetilde{\Phi}}
\newcommand{\sph}{\mathbb{S}}
\newcommand{\cerc}{\mathbb{S}^1}
\newcommand{\trho}{\widetilde{\rho}}
\newcommand{\tA}{\widetilde{A}}

\newcommand{\vreg}{V_{\reg}}
\newcommand{\breg}{B_{\reg}}
\newcommand{\vvois}{V_{\vois}}

\newcommand{\cu}{\mathcal{U}}
\newcommand{\cv}{\mathcal{V}}

\DeclareMathOperator{\can}{can}
\DeclareMathOperator{\std}{std}
\DeclareMathOperator{\Sp}{Sp}
\DeclareMathOperator{\ident}{id}
\DeclareMathOperator{\modulo}{mod}
\DeclareMathOperator{\im}{Im}
\DeclareMathOperator{\pair}{even}
\DeclareMathOperator{\impair}{odd}
\DeclareMathOperator{\Id}{Id}
\DeclareMathOperator{\reg}{reg}
\DeclareMathOperator{\vois}{neigh}
\DeclareMathOperator{\Ker}{ker}

\title{Families of monotone symplectic manifolds constructed via 
symplectic cut and their Lagrangian submanifolds} 
\author{Agn\`es GADBLED
\thanks {
\newline Partially supported by ANR project "Floer Power" 
ANR-08-BLAN-0291-03, SNF project "Complexity and recurrence in Hamiltonian 
systems" and MSRI post-doc fellowship NSF grant DMS-0441170}}
\date{January 2010}

\begin{document}
\selectlanguage{english}
\maketitle

\begin{abstract}
  We describe families of monotone symplectic manifolds 
constructed via the symplectic cutting procedure of Lerman~\cite{MR1338784}
from the cotangent bundle of manifolds endowed with a free circle action.
We also give obstructions to the monotone Lagrangian embedding 
of some compact manifolds in these symplectic manifolds.
\end{abstract}

\section*{Introduction}
\label{sec:intro}

We describe families of examples of monotone symplectic manifolds 
constructed by the symplectic cutting procedure of Lerman~\cite{MR1338784}
from the cotangent bundle of manifolds endowed with a free circle action.

If $V$ is a closed manifold endowed with a free circle action, this action
induces a Hamiltonian action of the circle on the cotangent bundle 
$\ctg V$ of $V$. If $H$ denotes the Hamiltonian associated to this action,
$W_\xi$ the symplectic manifold obtained by symplectic cut of $\ctg V$ at 
the level $\xi$ and $Q_\xi$ the reduced manifold~$H^{-1}(\xi)/\cerc$, 
we prove that $W_\xi$ has the following properties:
\begin{description}

\item[1)] (Lemma~\ref{lemfibr'eass}) The manifold $W_\xi$ is a complex line 
bundle over $Q_\xi$.

\item[2)] (Proposition~\ref{prop:omega-reduc}) The quotient $Q_\xi$ is 
a symplectic submanifold of~$W_\xi$ diffeomorphic to the cotangent 
bundle of $B=V/\cerc$. 
If $e$ is the Euler class of the circle bundle $V \rightarrow B$ and 
$q$ is the projection $Q_\xi \rightarrow B$, then the class of the 
symplectic form over $Q_\xi$ is $- 2 \pi \xi q^\star e$.

\item[3)] (Proposition~\ref{prop:mon-W}) The class in $H^2(W_\xi;\rr)$ 
of the symplectic form of~$W_\xi$ is related to the first Chern class 
of $W_\xi$ by the formula:
$$[\omega_\xi] = - 2 \pi \xi c_1(W_\xi).$$
In particular, $W_\xi$ is monotone for negative $\xi$.

\item[4)] (Proposition~\ref{prop:per}) The manifold $W_\xi$ is endowed 
with a Hamiltonian circle action and the sum of the weights of 
this action in a fixed point is equal to~$1$.

\item[5)] (Theorem~\ref{thm:Vmonotone}) The zero section $V$ of the 
cotangent bundle can be embedded as a monotone Lagrangian submanifold 
of $W_\xi$ for any negative level $\xi$.

\item[6)] (Corollary~\ref{cor:exactedsctgdsspcut}) Any exact compact 
Lagrangian submanifold of $\ctg V$ with vanishing Maslov class 
can be embedded as a monotone Lagrangian submanifold of $W_\xi$ 
for some negative level $\xi$.
\end{description}

In the general case of a monotone Lagrangian submanifold in a symplectic 
manifold, Seidel~\cite{MR1765826} has proved, 
under some assumptions on the symplectic manifold and the Lagrangian 
submanifold, that the Floer homology can be absolutely graded. He also 
proved that if the symplectic manifold is endowed with a Hamiltonian 
circle action, the Floer homology of the Lagrangian submanifold 
is periodic with respect to the absolute grading with period 
the sum of the weights of the action in a fixed point. 
Applying these results to the monotone Lagrangian submanifolds 
of the symplectic cut $W_\xi$, we get:
\begin{description}
\item[1)] (Theorem~\ref{thm:nsc-spcut}) An obstruction result on the 
embedding of a compact simply-connected manifold as a monotone Lagrangian 
submanifold in a symplectic cut under conditions on the Euler 
class $e$.

\item[2)] (Theorems \ref{thm:nsc} and \ref{thm:cassc}) 
Results on the index of the fundamental group of 
an exact compact Lagrangian submanifold with vanishing Maslov class in 
the cotangent bundle of a manifold endowed with a free circle action.

\item[3)] (Theorem~\ref{thm:tore-mondsspcut}) The monotone Lagrangian tori 
in a simply connected symplectic cut have a Maslov number equal to $2$.\\
\end{description}

In the first section, we explain the symplectic cutting construction
for our specific circle action and describe it in some examples. In
the second section, we study the structure of the spaces which appear 
in the construction and their monotonicity. The third part is 
essentially devoted to the monotonicity of the canonical embedding 
of~$V$ in 
the symplectic cut. After briefly recalling Seidel's periodicity results, 
we apply them to the monotone Lagrangian submanifolds of the 
symplectic cut in the fourth part and prove the results stated 
above. An appendix contains complements to some proofs and a 
discussion of the monotonicity in the case of a semi-free circle 
action.

\begin{acknowledgement}
\emph{This is a part of the author's PhD thesis at the University of Strasbourg. 
I would like to thank Mich\`ele Audin and Mihai Damian for their
  helpful suggestions and comments and for their careful reading of a earlier 
version of this article.
 I am greatly indebted to them for the time they spend discussing this
 subject with me.}
\end{acknowledgement}

\section{Symplectic cut of cotangent bundles}
\label{sec:SC-cotg}

Let $\pi : V \rightarrow B$ be a principal circle bundle, 
the base $B$ being path connected. Denote by
$X$ the fundamental vector field associated with
the circle action on~$V$. As this $\sph^1$-action is free, the
fundamental vector field cannot vanish.

\begin{rem}
We refer to appendix~\ref{sec:Append-Semilibre} for the study 
of the case of a semi-free circle action.
\end{rem}

Denote by $e \in H^2(B;\zz)$ the Euler class of the principal
bundle  $\pi$ and by $N_{e}$ the nonnegative
generator of the subgroup $\langle e , \pi_{2}(B) \rangle$ of $\zz$.\\

\textbf {Example 0:} In the case of the trivial bundle 
$V=B \times \sph^1 \rightarrow B$, $\pi$ is the projection on the first
factor, the circle action is the complex multiplication on the factor
$\sph^1$ (seen as the group of complex numbers of module $1$) and
the Euler class is zero.

A subcase is the case of $\sph^1$ acting on itself by multiplication, the
quotient $B$ being reduced to a point.\\

\textbf {Example 1:} If $V=\sph^{2n+1} \rightarrow B=\pnc$ ($n \geq 1$)
is the Hopf bundle, the Euler class is equal to the opposite of the
preferred generator of $H^2(\pnc;\zz)$ (this can also be considered as
a choice for the preferred generator of $H^2(\pnc;\zz)$).\\

\textbf {Example 2:} Let $V$ be the lens space $L_{p}^{2n+1}$,
quotient of the sphere $\sph^{2n+1}$ by the action of the subgroup
$\zz/p$ of $\sph^{1}$ consisting of the $p$-roots of unity.
$$
\xymatrix{
  \sph^{2n+1} \ar[rd] \ar[dd]\\ & L_{p}^{2n+1}  \ar[ld] \\
  \pnc& }
$$
The free action of $\sph^{1}$ on $\sph^{2n+1}$ enables to define a free
action of $\sph^{1}$ on $L_{p}^{2n+1}$, with quotient $\pnc$, such that the
Euler class of the principal circle bundle $L_{p}^{2n+1} \rightarrow 
\pnc$ is $-p$ times the preferred generator of $H^2(\pnc;\zz)$ (chosen in 
Example 1).\\

\textbf {Example 3:} We can also consider the Stiefel manifold 
$V_2(\rr^{n+2})$ consisting of the pairs of orthonormal
vectors of $\rr^{n+2}$. This manifold $V_2(\rr^{n+2})$ can also be
described as the quotient $SO(n+2)/SO(n)$. The subgroup $SO(2)
\simeq \sph^{1}$ of $SO(n+2)$ acts freely on $V_2(\rr^{n+2})$, the
quotient of this action being the Grassmannian $B =
\widetilde{G}_2(\rr^{n+2})$ of the ($2$-dimensional) oriented planes
in $\rr^{n+2}$.

This Grassmannian can be identified with the quadric $Q^{n}$ in
$\mathbb{C} \mathbb{P}^{n+1}$ of equation 
$\displaystyle \sum_{j=0}^{n+1} z_{i}^{2}=0$. To see this, we can 
use the first description of the Stiefel manifold $V_2(\rr^{n+2})$ as
a subset of $(\rr^{n+2})^{2}$:
$$V_2(\rr^{n+2}) = \{(x,y) \in (\rr^{n+2})^{2} |\; x \cdot y = 0, \|x\|=\|y\|=1 \}.$$
The action of $SO(2) \simeq \sph^1$ can then be written: 
$$e^{i \theta} \cdot (x,y)=(\cos(\theta) x-\sin(\theta) y,\sin(\theta) 
x+\cos(\theta) y).$$
The map 
$$
\begin{array}{ccc}
V_2(\rr^{n+2}) & \longrightarrow & \mathbb{C} \mathbb{P}^{n+1} \\
      (x,y)    & \longmapsto     & [x+iy]
\end{array}
$$
descends to the quotient and defines an embedding of the Grassmannian 
$\widetilde{G}_2(\rr^{n+2})$ in $\mathbb{C} \mathbb{P}^{n+1}$, 
the image of which is the quadric $Q^{n}$ (see \cite{MR2357797}). 

For $n \geq 3$, by the Lefschetz hyperplane theorem, the
embedding of the quadric induces an isomorphism between 
$H^{2}(\widetilde{G}_2(\rr^{n+2});\zz)$ and  
$H^{2}(\mathbb{C} \mathbb{P}^{n+1};\zz)$. Therefore, 
$H^{2}(\widetilde{G}_2(\rr^{n+2});\zz)$ is isomorphic to $\zz$, 
a preferred generator being given by the pull-back of the preferred 
generator of  
$H^{2}(\mathbb{C} \mathbb{P}^{n+1};\zz)$. Besides, $V_2(\rr^{n+2})$ is
the total space of the pullback of the Hopf bundle $\sph^{2n+3}
\rightarrow \mathbb{C} \mathbb{P}^{n+1}$ by the embedding $j$ of
$\widetilde{G}_2(\rr^{n+2})$ in $\mathbb{C} \mathbb{P}^{n+1}$:
$$
\xymatrix{
j^{\star}\sph^{2n+3}=V_2(\rr^{n+2}) \ar[r] \ar[d] & \sph^{2n+3} \ar[d] \\
\widetilde{G}_2(\rr^{n+2})   \ar[r]^-{j}    & \mathbb{C} \mathbb{P}^{n+1}
           }
$$
Indeed, with the map $V_2(\rr^{n+2}) \rightarrow \sph^{2n+3}$ given by 
$\displaystyle (x,y) \mapsto \frac{1}{\sqrt{2}}(x+iy)$ and the
projection  $V_2(\rr^{n+2}) \rightarrow \widetilde{G}_2(\rr^{n+2})$,
the above diagram commutes. We deduce by naturality of the Euler class
that:
\begin{eqnarray*}
  e(V_2(\rr^{n+2}) \rightarrow \widetilde{G}_2(\rr^{n+2})) 
  & = & j^{\star}e(\sph^{2n+3} \rightarrow \mathbb{C} \mathbb{P}^{n+1})
\end{eqnarray*}
and hence, for $n \geq 3$, the Euler class of the bundle
$V_2(\rr^{n+2}) \rightarrow \widetilde{G}_2(\rr^{n+2})$ is equal to
the opposite of the generator of $\widetilde{G}_2(\rr^{n+2})$.\\
For the low dimensions $n=1$ and $n=2$, the description of the Stiefel 
manifolds and the corresponding Grassmannians are well known (see for 
instance~\cite{MR759162}). The long exact sequence associated with the
principal circle bundle 
$V_2(\rr^{n+2}) \rightarrow \widetilde{G}_2(\rr^{n+2})$ gives the
numbers: for $n=1$, $N_{e}=2$ and for $n=2$, $N_{e}=1$.\\

\subsubsection*{Associated Hamiltonian action on the cotangent bundle}

For $u \in \sph^{1}$, denote by $\rho(u): V \rightarrow V$ (defined by
$\rho(u)(x) = u \cdot x$) the action of~$\sph^{1}$ on $V$. With this
circle action on $V$ is associated a $\sph^1$-action on $\cotg$:
$$
\begin{array}{ccc}
\sph^1 \times \cotg & \longrightarrow & \cotg \\
  (u,(x,\varphi))& \longmapsto     & (u \cdot x, u \cdot \varphi)
\end{array}
$$
where  
$$u \cdot \varphi = \,^{t}(T_x \rho(u)^{-1})(\varphi),$$
namely, if $v \in T_{u \cdot x} V$, then 
$$u \cdot \varphi = \varphi (T_x \rho(u)^{-1} v).$$
If $\cotg$ is endowed with its canonical symplectic structure, this
action is Hamiltonian with Hamiltonian function:
$$
\begin{array}{cccc}
H : &   \cotg     & \longrightarrow & \rr \\
    & (x,\varphi) & \longmapsto     & \langle \varphi , X(x) \rangle
\end{array}
$$
(see for instance \cite [exercise 3.12]{MR1373431}).

This action is free since the circle action is free on $V$ and
the canonical projection $p_V : \cotg \rightarrow V$ is an equivariant
map.

In particular, this action has no fixed point, the associated
Hamiltonian vector field $X_H$ does not vanish
(because the projection of $X_H$ on $TV$ by $T p_V$ is $X$, a vector
field that never vanishes) and all the values of the Hamiltonian are
regular.

\subsubsection*{Symplectic cut of the cotangent bundle}

A symplectic manifold endowed with a Hamiltonian circle action being given, 
Lerman~\cite{MR1338784} has developped a construction which enables to embed 
the reduced spaces as codimension $2$ symplectic submanifolds of a symplectic 
manifold (called symplectic cut) 

We describe this construction in the case of the cotangent bundle $\cotg$. This 
symplectic cut is naturally endowed with a Hamiltonian circle action 
(Proposition~\ref{prop:per}). We are interested in the monotonicity of this 
symplectic manifold (Section~\ref{sec:Struct-des-varits}) in order to apply 
a result of Seidel~\cite{MR1765826} (reminded in Section~\ref{sec:Seidel}) on 
the periodicity of the Lagrangian Floer homology.

To define the symplectic cut construction, we consider the following 
Hamiltonian action of $\sph^1$ on the product $\cotg \times \cc$ : 
$$
\begin{array}{ccc}
\sph^1 \times (\cotg \times \cc) & \longrightarrow & \cotg \times \cc \\
    (u,(x,\varphi,z))         & \longmapsto     & (u \cdot x, u \cdot \varphi, \bar{u} z)
\end{array}
$$
with associated Hamiltonian:
$$
\begin{array}{cccc}
\tH : & \cotg \times \cc & \longrightarrow & \rr \\
      &  (x,\varphi,z)  & \longmapsto     & H(x,\varphi) - \frac{1}{2} | z |^2.
\end{array}
$$
As for $H$, this action is free and the levels of the Hamiltonian $\tH$
are regular.

Let $\xi \in \rr$ be a (regular) value of $\tH$. The level
$\tH^{-1}(\xi)$ is the disjoint union of 
$$\{(x,\varphi,0) |H(x,\varphi) = \xi \},$$
diffeomorphic to $H^{-1}(\xi)$, and  
$$\{(x,\varphi,z) | H(x,\varphi) > \xi \mbox{ and } |z|=\sqrt{2(H(x,\varphi) - \xi)}\},$$ 
diffeomorphic to $H^{-1}(\xi, + \infty) \times \sph^1$, the
$\sph^{1}$-action preserving both sets.

Notice that $\sph^1$ acts also freely on the level $H^{-1}(\xi)$ since 
$$
\begin{array}{cccc}
\widetilde{q}_{\xi}: & H^{-1}(\xi) & \longrightarrow & V \\
                 & (x,\varphi) & \longmapsto     & x
\end{array}
$$
is an equivariant fibration (see Lemma~\ref{lemtq}) and the
action is free on the base space.

\begin{figure}[htbp]
  \begin{center}
   \psfrag{Vg}{$V$}
   \psfrag{T*V}{$\cotg$}
   \psfrag{xi}{$\xi$}
   \psfrag{0}{$0$}
   \psfrag{H}{$H$}
   \psfrag{Vd}{$V$}
   \psfrag{Wxi}{$W_{\xi}$}
   \psfrag{Qxi}{$Q_{\xi}$}
   \includegraphics{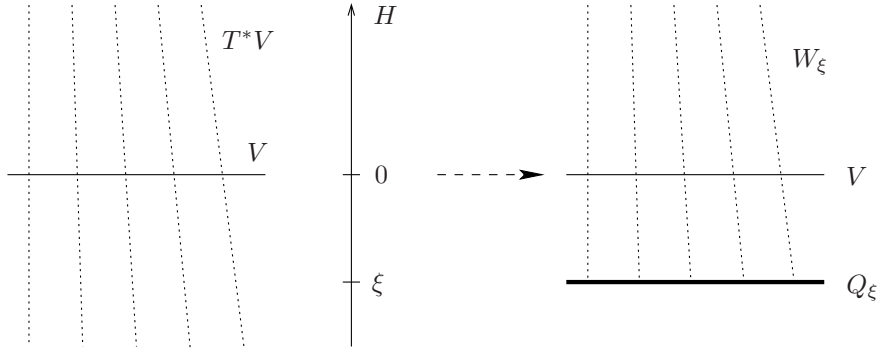}
    \caption{The cotangent bundle (left) and its symplectic cut
      (right)} 
    \label{fig:cotangent}
  \end{center}
\end{figure}

As the action of $\sph^1$ on the regular level $\tH^{-1}(\xi)$ is free,
we can carry out a symplectic reduction. Equipped with the reduced
symplectic form, the quotient $W_{\xi} = \tH^{-1}(\xi) / \sph^1$ is a
symplectic manifold. The decomposition above descends to the
quotient, the image of the part $H^{-1}(\xi)$ is a symplectic 
submanifold $Q_{\xi}$ of $W_{\xi}$ diffeomorphic to
$H^{-1}(\xi) /\sph^1$. Its complement is  symplectically diffeomorphic to
the open set $H^{-1}(\xi, +\infty)$ (see Figure
\ref{fig:cotangent}). 

\begin{rem}
\label{remxicot}
If $\xi < 0$, then $H^{-1}(\xi, + \infty)$ is a neighbourhood of the
zero section in the cotangent bundle $\cotg$ and this neighbourhood is 
symplectically embedded as an open submanifold of $W_{\xi}$. In
particular, $V$ is canonically embedded as a Lagrangian submanifold 
of the symplectic cut $W_{\xi}$.
\end{rem}

Because of Remark~\ref{remxicot} and of the monotonicity condition on 
$W_{\xi}$ (Proposition~\ref{prop:mon-W}), we will deal in the following 
only with the case $\xi < 0$.\\

\textbf {Example 0:} Let us first study the case of the circle 
$V=\sph^1$ acting on itself by multiplication with quotient $B$ reduced 
to a point. The fundamental field of the action at the point $x$
is then $ix$ (we consider $\sph^{1}$ as the unit circle in $\cc$). 
The cotangent space
$T^{\ast}\sph^{1}$ is diffeomorphic to $\sph^{1} \times \rr$: if $x$ is a point
of $\sph^{1}$ and $\varphi$ is an element of $T_{x}^{\ast}\sph^{1}$, there
exists a unique real number $\lambda$ such that 
$\varphi= \langle \lambda \; iu,\cdot \rangle$ (denoting by 
$\langle \cdot,\cdot \rangle$ the scalar product on $\cc=\rr^{2}$). An 
explicit diffeomorphism from $T^{\ast}\sph^{1}$ to $\sph^{1} \times \rr$ can
be given by mapping the pair $(x,\varphi)$ on the pair $(x,\lambda)$.

The circle action on $T^{\ast}\sph^{1}$ can then be written on 
$\sph^{1} \times \rr$:
$$u \cdot (x,\lambda) = (ux, \lambda)$$
and the Hamiltonians become:
$$H(x,\lambda) = H(x,\varphi)=\varphi(ix)=\lambda$$
and 
$$\tH(x,\lambda,z) = H(x,\varphi,z)=\lambda-\frac{1}{2}|z|^{2}.$$
Thus, the level $H^{-1}(\xi) = \{(x,\xi)| \: x \in \sph^{1}\}$ is a circle
and the quotient $Q_{\xi}=H^{-1}(\xi)/\sph^{1}$ is a point.

For $\tH$, the level $\tH^{-1}(\xi)$ is the set  
$\{ ( x,\xi+\frac{1}{2}|z|^{2},z ) | \: x \in \sph^{1}, z \in \cc \}$ 
and the quotient $W_{\xi}= \tH^{-1}(\xi)/\sph^{1}$ is symplectomorphic to
$\cc$. Indeed, the one-to-one map 
$$
\begin{array}{cccc}
i_{\cc} : & \cc & \longrightarrow & W_{\xi} \\
          &  z  & \longmapsto     & \left[{ 1,\xi+\frac{1}{2}|z|^{2},z }\right]
\end{array}
$$
is also onto as each element 
$$\left[{x,\xi+\frac{1}{2}|z|^{2},z}\right]=\left[{1,\xi+\frac{1}{2}|z|^{2},x^{-1}z}\right]$$
of $W_{\xi}$ can be written as the image by $i_{\cc}$ of a complex
number ($x^{-1}z=\bar{x}z$ here).

Moreover, the symplectic form $\omega_{\xi}^{W}$ on $W_{\xi}$ is the
symplectic reduction of the $2$-form 
$\omega_{\can} \oplus \omega_{\std}$ where 
$\omega_{\can}$ is the canonical symplectic form on $T^{\ast}\sph^{1}$
and $\omega_{\std}$ is the standard symplectic form on $\cc$. But if
$\widetilde{\imath}$ is the embedding of $\cc$ in $\tH^{-1}(\xi)$ defined as : 
$$\widetilde{\imath}(z)=\left( {1,\xi+\frac{1}{2}|z|^{2},z} \right),$$
then 
\begin{eqnarray*}
(i^{\star}\omega_{\xi}^{W})_{z}(\zeta_{1},\zeta_{2}) 
& = & (\, \widetilde{\imath}^{\, \star}(\omega_{\can} \oplus
\omega_{\std}))_{z}(\zeta_{1},\zeta_{2})\\
& = & \omega_{\std}(\zeta_{1},\zeta_{2})
\end{eqnarray*}
since $z \mapsto (1,\xi+\frac{1}{2}|z|^{2})$ takes value in a fibre of
the projection $T^{\ast}\sph^{1} \rightarrow \sph^{1}$ and the fibres are
Lagrangian submanifolds.

\begin{figure}[htbp]
  \centering
   \psfrag{V}{$V$}
   \psfrag{H}{$H$}
   \psfrag{R}{$\rr$}
   \psfrag{xi}{$\xi$}
   \psfrag{0}{$0$}
   \psfrag{Wxi}{$W_{\xi}=\cc$}
   \psfrag{Qxi}{$Q_{\xi}$}
   \includegraphics{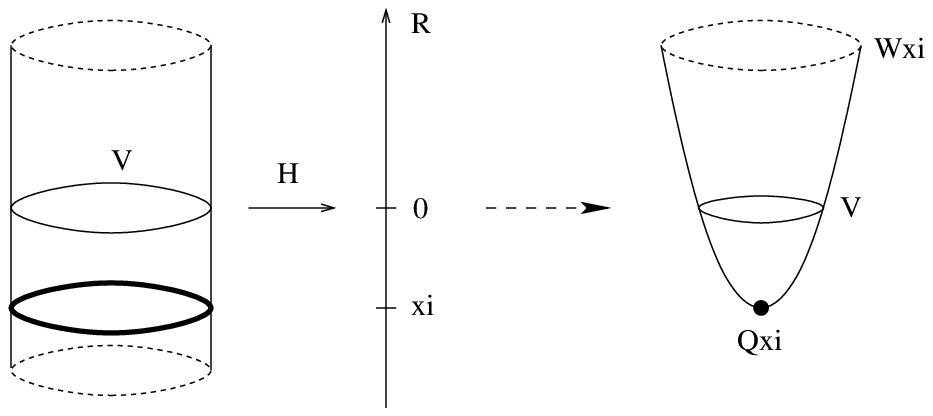}
  \caption{Symplectic cut of $T^\ast \sph^1$}
  \label{fig:SCcerc}
\end{figure}

The manifold $V = \sph^{1}$ is embedded as the zero section in $\cotg$ 
and as the image of the map 
$$ x \mapsto [x,0,\sqrt{-2 \xi}] = [1,0,\bar{x} \sqrt{-2
  \xi}]$$
in $W_{\xi}$ if $\xi < 0$. Its image in $\cc$ by the symplectic
diffeomorphism just described is the circle centered at the origin and 
of radius $\sqrt{-2 \xi}$. This circle bounds a disc of area 
$-2 \pi \xi$ (which is also the area of the cylinder that lies between
the zero section and the level $\xi$) (see Figure~\ref{fig:SCcerc}).

In the general case of a product $V=B \times \sph^1 \rightarrow B$, there
is a symplectic diffeomorphism $\cotg \simeq T^{\ast}B \times
T^{\ast}\sph^{1}$. The circle action on $V$ being concentrated on the
factor $\sph^1$, 
the action on $\cotg$ is concentrated on $T^{\ast}\sph^{1}$ and the
corresponding Hamiltonian on $T^{\ast}B \times T^{\ast}\sph^{1}$ can be
written, for any $((b,\varphi),(x,\lambda)) \in T^{\ast}B \times
T^{\ast}\sph^{1}\simeq T^{\ast}B \times (\sph^{1} \times \rr),$
$$H((b,\varphi),(x,\lambda)) = \lambda.$$
Thus, the symplectic cut of $\cotg$ (at the level $\xi$) is the
product of $T^{\ast}B$ and of the symplectic cut of $T^{\ast}\sph^{1}$
(also at the level $\xi$) and hence is symplectomorphic to 
$T^{\ast}B \times \cc$.

If $\xi$ is negative, $V = B \times \sph^1$ can be embedded as the
product of the zero section of $T^{\ast}B$ and of the circle centered
at the origin of radius $\sqrt{-2 \xi}$.\\

Let us now analyse the structure of these manifolds. In particular, we
would like to know if $W_{\xi}$ is monotone
(see Section \ref{sec:Struct-des-varits}) and in this case, if the
zero section is monotone in $W_{\xi}$ (Section
\ref{sec:MonotonieV}).\\

\section{Monotonicity of the symplectic cut}
\label{sec:Struct-des-varits}

Recall (see \cite{MR1223659} and \cite{MR1367384}) that a symplectic
manifold $(W,\omega)$ is said to be monotone if there exists a positive 
constant $K_W$ such that for any sphere $w: \sph^2 \rightarrow W$, 
\begin{equation}
\int w^{\star}\omega = K_W \, c_1([w]),
\label{eq:W-mon}
\end{equation}
where $c_1$ is the first Chern class of the tangent bundle $TW \rightarrow W$. 

If $W$ is a symplectic manifold, denote by $N_W$ its first Chern
number, namely the nonnegative generator of the subgroup 
$\langle c_1(W), \pi_2(W) \rangle$ of $\zz$.\\

We will prove that the symplectic structure and the first Chern class
of $W_{\xi}$ can be deduced of those of $Q_{\xi}$. To begin with, let
us investigate the reduced manifold $Q_{\xi} = H^{-1}(\xi) /\sph^1$.\\

\subsection{The reduced symplectic manifold $Q_{\xi}$}
\label{sec:etude-q-xi}

\begin{lem}
\label{lemtq} 
The projection $\widetilde{q}_{\xi} : H^{-1}(\xi) \rightarrow V$ is an
equivariant fibration and descends to the quotients in a fibration
$$
\begin{array}{cccc}
q_{\xi}: & Q_{\xi}     & \longrightarrow & B \\
         & [x,\varphi] & \longmapsto     & \pi(x).
\end{array}
$$
\end{lem}

\begin{proof}
Let us begin by looking at this projection for $\xi = 0$. In this
case, 
$$
\begin{array}{cccc}
\widetilde{q}_0: & \{(x,\varphi) | \langle \varphi , X(x) \rangle = 0\}   & \longrightarrow & V \\
                  & (x,\varphi)                                            & \longmapsto     & x
\end{array}
$$
is a vector subbundle of the cotangent bundle $\cotg \rightarrow V$: it is
the annihilator of the vector field $X$ in the cotangent
bundle. Descending to the quotients, $q_0$ is also a 
vector bundle over $B$.

Let $\{\alpha_x\}_{x \in V}$ be a family of linear forms $\alpha_x \in
T_{x}^{\ast}V$ such that for all $x \in V$, $\langle \alpha_x , X(x)
\rangle = 1$ and for all $u \in \sph^1$, $\alpha_{u \cdot x} = u \cdot
\alpha_x$. For instance, it can be defined using an $\sph^{1}$-invariant
Riemannian metric $g$ on $V$ and setting 
$\alpha_x = g_{x}( \, \cdot \, , X(x) )$.
For any $\xi$, this allows to define an equivariant diffeomorphism 
$$
\begin{array}{cccc}
\tPhi_{\xi}: & H^{-1}(\xi) & \longrightarrow & H^{-1}(0) \\
            & (x,\varphi) & \longmapsto     & (x,\varphi - \xi \alpha_x).
\end{array}
$$
such that $\widetilde{q}_{\xi} = \widetilde{q}_0 \circ \tPhi_{\xi}$.

Descending to the quotients, the diffeomorphism $\tPhi_{\xi}$ defines a
diffeomorphism $\Phi_{\xi} : Q_{\xi} \rightarrow Q_0$ such that 
$q_{\xi} = q_0 \circ \Phi_{\xi}$. \end{proof}

\begin{cor}
The level $H^{-1}(0)$ (and hence $H^{-1}(\xi)$ as well) has the
homotopy type of $V$ and the quotient $Q_0$ (and hence $Q_{\xi}$ as
well) has the homotopy type of~$B$, the homotopy equivalence being
given by the projections. 
\end{cor}

Furthermore, the fibre bundle $q_0$ is a particular fibre bundle on
$B$:

\begin{lem}
\label{lem:encyclo}
The space $Q_0$ endowed with the symplectic form obtained by reduction
of the canonical symplectic form on $\cotg$ is symplectomorphic to the
cotangent bundle $T^{\ast}B$.
\end{lem}
 
A proof of Lemma \ref{lem:encyclo} can be found in \cite[Section 5]{magnetic}.

\begin{cor}
\label{cor:Chern-Qxi}
The first Chern class $c_1(TQ_{\xi} \rightarrow Q_{\xi})$ of $Q_{\xi}$
is zero in $H^2(Q_{\xi};\rr)$.
\end{cor}

\begin{proof}
  Thanks to Lemma~\ref{lem:encyclo}, $Q_0$ is symplectomorphic to~$T^{\ast}B$.
  As the first Chern class 
  $c_1 \left({T(T^{\ast}B) \rightarrow T^{\ast}B }\right)$ is 
  a $2$-torsion element in $H^2(T^{\ast}B;\zz)$ (because the restriction 
  to $B$ of the tangent bundle of $\ctg B$ is isomorphic to the 
  complexified bundle of a real bundle, see~\cite{MR0440554}), 
  $c_1(TQ_{0} \rightarrow Q_{0})$ is $2$-torsion element 
  in $H^2(Q_{0};\zz)$.
 
  Moreover, $Q_{\xi}$ varies continuously with $\xi$ and the first
  Chern class of $Q_{\xi}$ takes values in a discrete space so that
  $c_1(TQ_{\xi} \rightarrow Q_{\xi})$ is also a $2$-torsion element 
  of~$H^2(Q_{\xi};\zz)$ and hence vanishes in $H^2(Q_{\xi};\rr)$.
\end{proof}

We determine the class $[\omega_{\xi}]$ of the
symplectic form on $Q_{\xi}$ with the Duistermaat-Heckman formula (see
\cite{MR674406}, we follow here the sign conventions of
\cite{MR2091310}):

\begin{prop}
\label{prop:omega-reduc}
The class in $H^2(Q_{\xi_0};\rr)$ of the reduced symplectic form 
on~$Q_{\xi}$ is equal to:
$$[\omega_{\xi}] =  - 2 \pi \xi \; q_{\xi}^{\star}e.$$
\end{prop}

\begin{proof}
As all the levels are regular for $H$, by the
Duistermaat-Heckman formula, for all $\xi$, $\xi_0$ in $\rr$, 
$$[\omega_{\xi}] - [\omega_{\xi_0}] = (\xi - \xi_0) \left( - 2 \pi e ( H^{-1}(\xi) \rightarrow Q_{\xi} ) \right)$$
in $H^2(Q_{\xi_0};\rr)$ where $e ( H^{-1}(\xi) \rightarrow Q_{\xi} )$
is the Euler class of the principal circle bundle: 
$\pi_Q: H^{-1}(\xi) \rightarrow Q_{\xi}$.

We have now to compute this Euler class 
$e ( H^{-1}(\xi) \rightarrow Q_{\xi} )$. As the bundle
$H^{-1}(\xi) \rightarrow Q_{\xi}$ is the pullback of the fibre bundle
$\pi: V \rightarrow B$ by the projection $q_{\xi}: Q_{\xi} \rightarrow
B$: 
$$
\xymatrix{
H^{-1}(\xi) \ar[r]^-{\widetilde{q}_{\xi}} \ar[d]^-{\pi_{Q}} & V \ar[d]^-{\pi} \\
Q_{\xi} \ar[r]^-{q_{\xi}}                                   & B
         }
$$
we have by naturality:
$$e (H^{-1}(\xi) \rightarrow Q_{\xi}) = q_{\xi}^{\star} \left( e (V \rightarrow B) \right) = q_{\xi}^{\star}e,$$ 
where $q_{\xi}^{\star}: H^2(B;\zz) \rightarrow H^2(Q_{\xi};\zz)$ is an 
isomorphism.

Eventually, the Duistermaat-Heckman formula gives:
$$[\omega_{\xi}] - [\omega_{\xi_0}] = (\xi - \xi_0) (- 2 \pi q_{\xi}^{\star}e).$$
For $\xi_0 = 0$, $(Q_0,\omega_{0})$ is symplectomorphic to the
cotangent bundle $(T^{\ast}B, \omega_{\can})$. But $\omega_{\can}$ is
an exact form, hence the class $[\omega_{\can}]$ and a fortiori the
class $[\omega_{0}]$ vanish, so that the Duistermaat-Heckman formula 
becomes with $\xi_0 = 0$:
$$[\omega_{\xi}] =  - 2 \pi \xi \; q_{\xi}^{\star}e.$$
\end{proof}

We could have also used \cite[section 5]{magnetic} and describe $Q_{\xi}$ 
as a magnetic cotangent bundle. Recall that a magnetic cotangent bundle is 
a cotangent bundle $q: \ctg B \rightarrow B$ whose symplectic form differs 
from the canonical symplectic form by a perturbation $q^\star \sigma$ where 
$\sigma$ is a closed $2$-form on $B$. Here $\sigma = \xi F(\alpha)$ where 
$F(\alpha)$ is the curvature of a connexion $1$-form $\alpha$ on the fibre 
bundle $\pi : V \rightarrow B$.\\

\textbf {Example 0:} In the case of the trivial bundle, $Q_{\xi}$ is
symplectomorphic to~$T^{\ast}B$ for all $\xi$, as the component in
$T^{\ast}B$ of the $1$-form $\alpha$ vanishes. The $2$-form $\omega_{\xi}$
is always exact and $[\omega_{\xi}] = 0$ in $H^2(Q_{\xi_0};\rr)$ as
confirmed by Proposition~\ref{prop:omega-reduc}.\\

\textbf {Examples 1, 2 and 3:} In the case of the Hopf bundle for
$n \geq 1$ (and similarly in the case of the Stiefel manifolds for $n
\geq 3$), the class $[\omega_{\xi}]$ is equal to $-2 \pi$ times the pullback on
$T^{\ast}\pnc$ of the image in $H^{2}(\pnc;\rr)$ of the 
generator of $H^{2}(\pnc;\zz)$ (respectively the pullback on
$T^{\ast}\widetilde{G}_2(\rr^{n+2})$ of the image of the 
generator of $H^{2}(\widetilde{G}_2(\rr^{n+2});\zz)$). In the case of
example 2, this class is $-2 \pi p$ times the generator.\\

\subsection{Consequences for the symplectic cut}

In this paragraph, we want to prove the following:

\begin{prop} 
\label{prop:mon-W}
If $V \rightarrow B$ is a principal circle bundle and if $\xi < 0$,
then the symplectic cut $W_{\xi}$ is monotone and contains $V$ as a
Lagrangian submanifold.
\end{prop}

This proposition is a consequence of the bundle structure of
$W_{\xi}$:

\begin{lem} 
\label{lemfibr'eass}
The quotient $W_{\xi}$ is the total space of a complex line bundle 
on~$Q_{\xi}$, namely the complex line bundle associated with the 
principal circle bundle: $H^{-1}(\xi) \rightarrow Q_{\xi}$.
\end{lem}

\begin{proof}
By definition, the complex line bundle associated with the principal
circle bundle $H^{-1}(\xi) \rightarrow Q_{\xi}$ is 
$H^{-1}(\xi) \times_{\sph^1} \cc \rightarrow Q_{\xi}$ where $\sph^1$ acts on
$H^{-1}(\xi) \times \cc$ by:
$$u \cdot (x, \varphi, z) = (u \cdot x, u \cdot \varphi, \bar{u}z)$$
and the map 
$$
\begin{array}{ccc}
\tH^{-1}(\xi) & \longrightarrow & H^{-1}(\xi) \times \cc \\
(x,\varphi,z) & \longmapsto     & (x,\varphi - \frac{1}{2} | z |^2 \alpha_x, z)
\end{array}
$$
is an equivariant diffeomorphism for the actions of $\sph^1$. It descends 
to the quotients as a diffeomorphism 
$W_{\xi} \rightarrow H^{-1}(\xi) \times_{\sph^1} \cc$. 
\end{proof}

\begin{cor}
\label{cor:struct-de-w_xi}
\begin{description}
\item[(i)] The projection $p_{\xi}: W_{\xi} \rightarrow Q_{\xi}$ is a 
homotopy equivalence. As a consequence, it is also the case for  
$q_{\xi} \circ p_{\xi} : W_{\xi} \rightarrow B$.
\item[(ii)] The first Chern class of the bundle 
  $p_{\xi}: W_{\xi} \rightarrow Q_{\xi}$ is equal to the Euler
  class of the circle bundle: $H^{-1}(\xi) \rightarrow Q_{\xi}$. 
\item[(iii)] The bundle $p_{\xi}: W_{\xi} \rightarrow Q_{\xi}$ is the normal
  bundle of $Q_{\xi}$ in $W_{\xi}$. 
\end{description}
\end{cor}

\textbf {Proof of (iii): }
By Lemma \ref{lemfibr'eass}, $W_{\xi} \rightarrow Q_{\xi}$
is isomorphic to the bundle $H^{-1}(\xi) \times_{\sph^1} \cc \rightarrow
Q_{\xi}$. But (see \cite[remark 1.7]{MR1338784}), this bundle
$H^{-1}(\xi) \times_{\sph^1} \cc \rightarrow Q_{\xi}$ is the normal
bundle of the reduced space $Q_{\xi}$ in $W_{\xi}$.

Indeed, we have described a diffeomorphism 
$H^{-1}(\xi) \times \cc \stackrel{\sim}{\rightarrow} \tH^{-1}(\xi)$ 
that induces an isomorphism between the tangent bundles:
$$T \left( H^{-1}(\xi) \times \cc \right) \stackrel{\sim}{\rightarrow} T \tH^{-1}(\xi).$$ 
Hence, if $j_{\xi} : H^{-1}(\xi) \rightarrow \tH^{-1}(\xi)$
denotes the injection of $H^{-1}(\xi)$ in $\tH^{-1}(\xi)$, then 
$j_{\xi}^{\star} T \left(H^{-1}(\xi) \times \cc \right)
\stackrel{\sim}{\rightarrow} j_{\xi}^{\star} T \tH^{-1}(\xi)$.

But $j_{\xi}^{\star} T \left( H^{-1}(\xi) \times \cc \right)$ is
isomorphic as a bundle on $H^{-1}(\xi)$ to the Whitney sum of $T
H^{-1}(\xi)$  and of the trivial bundle $H^{-1}(\xi) \times \cc$.

All the maps above are equivariant under the circle action and
induce an isomorphism of bundles on $Q_{\xi}$:
$$i_{\xi}^{\star} T W_{\xi} \simeq T Q_{\xi} \oplus \left( H^{-1}(\xi) \times_{\sph^1} \cc \right),$$
denoting by $i_{\xi}: Q_{\xi} \rightarrow W_{\xi}$ the embedding of
$Q_{\xi}$ as the zero section of $W_{\xi}$.\hfill \qedsymbol\\

From Corollary \ref{cor:struct-de-w_xi} we get the first
Chern class of $T W_{\xi} \rightarrow W_{\xi}$:

\begin{lem}
The first Chern class of $W_{\xi}$ in $H^2(W_{\xi};\rr)$ is 
$$c_{1}(W_{\xi}) = (q_{\xi} \circ p_{\xi})^{\star} e.$$
\end{lem}

\begin{proof}
By part (iii) of Corollary \ref{cor:struct-de-w_xi}, 
$$c_1(i_{\xi}^{\star} T W_{\xi} \rightarrow Q_{\xi}) = i_{\xi}^{\star}
c_1(T W_{\xi} \rightarrow W_{\xi}) = c_1(T Q_{\xi} \rightarrow
Q_{\xi}) + c_1(W_{\xi} \rightarrow Q_{\xi})$$
and if we apply the inverse isomorphism $p_{\xi}^{\star}$ of
$i_{\xi}^{\star}$:
$$c_1(T W_{\xi} \rightarrow W_{\xi}) = p_{\xi}^{\star} \left({ c_1(T
    Q_{\xi} \rightarrow Q_{\xi}) + c_1(W_{\xi} \rightarrow Q_{\xi})
    }\right).$$
But, by Corollary \ref{cor:Chern-Qxi}, the first Chern class 
$c_1(T Q_{\xi} \rightarrow Q_{\xi})$ vanishes in $H^2(Q_{\xi};\rr)$
and the relation becomes in $H^2(W_{\xi};\rr)$:
\begin{eqnarray*}
c_1(TW_{\xi} \rightarrow W_{\xi}) &=& p_{\xi}^{\star} (c_1(W_{\xi} \rightarrow Q_{\xi}) \\
                                  &=& (q_{\xi} \circ p_{\xi})^{\star} e.
\end{eqnarray*}
\end{proof}

\begin{cor}
  The first Chern number $N_{W_{\xi}}$ of $W_{\xi}$ is equal to
  $N_{e}$ (the nonnegative generator of   $\langle e , \pi_{2}(B) \rangle$).
\end{cor}

The bundle structure of $W_{\xi} \rightarrow Q_{\xi}$ allows also
to determine the class of the symplectic form on $W_{\xi}$:

\begin{lem}
If $\omega^W_{\xi}$ denotes the symplectic form on $W_{\xi}$
defined by the symplectic reduction, then 
$$[\omega^W_{\xi}] = -2 \pi  \xi (q_{\xi} \circ p_{\xi})^{\star} e$$
in $H^2(W_{\xi};\rr)$.
\end{lem}

\begin{proof}
As $Q_{\xi}$ is a symplectic submanifold of $W_{\xi}$, we have
$$\omega_{\xi} = i_{\xi}^{\star} \omega^W_{\xi}$$ 
(still denoting by
$i_{\xi}$ the embedding of $Q_{\xi}$ in $W_{\xi}$) and hence 
$$[\omega_{\xi}] = i_{\xi}^{\star} [\omega^W_{\xi}] \mbox{ in }
H^2(Q_{\xi};\rr).$$ 
Applying the inverse isomorphism $p_{\xi}^{\star}$
of $i_{\xi}^{\star}$, we get 
$$[\omega^W_{\xi}] = p_{\xi}^{\star} [\omega_{\xi}] \mbox{ in } 
H^2(W_{\xi};\rr).$$ 
But we have proved that $[\omega_{\xi}]= - 2 \pi \xi
q_{\xi}^{\star}e$, thus
$$[\omega^W_{\xi}] = - 2 \pi \xi (q_{\xi} \circ p_{\xi})^{\star} e.$$
\end{proof}

\textit{Proof of Proposition \ref{prop:mon-W}.} 
We can now compare the class of the symplectic form $[\omega^W_{\xi}]$
and the first Chern class $c_1(W_{\xi})$ in $H^2(W_{\xi};\rr)$:
$$ [\omega^W_{\xi}] = - 2 \pi \xi \, c_1(W_{\xi})$$
Together with Remark~\ref{remxicot}, this equality
gives the proposition. \hfill \qedsymbol

\begin{rem}
It is also possible to carry out the symplectic cutting with the
Hamiltonian 
$$(x,\varphi,z) \mapsto H(x,\varphi) + \frac{1}{2} | z |^2$$ 
instead of
$\widetilde{H}$. In this case, it is the open subset 
$H^{-1}(-\infty, \xi)$ that embeds in the new symplectic cut
$W'_{\xi}$ and the first Chern class is 
$c_1(W'_{\xi} \rightarrow Q_{\xi}) = (q_{\xi} \circ p_{\xi})^{\star} (-e)$. 
Thus, it is for $\xi > 0$ that $W'_{\xi}$ is monotone and contains $V$
as Lagrangian submanifold. 
\end{rem}

Let us look again at our examples.\\

\textbf {Example 0:} In the case of the trivial bundle, the Euler
class and the class of the symplectic form are trivial in
$H^2(W_{\xi};\rr)$. We already know it, since $W_{\xi}$ is then
symplectomorphic to the product of a cotangent bundle and $\cc$.\\

\textbf {Example 1:} In the case of the Hopf bundle, the reduced
manifold $W_{\xi}$ is simply connected, with first Chern number
$N_{W_{\xi}}=1$. For $\xi < 0$, it is monotone and the sphere
$\sph^{2n+1}$ can be embedded in $W_{\xi}$ as a Lagrangian submanifold.\\ 

\textbf {Example 2:} In the case when $V = L^{2n+1}_p$ is the
$2n+1$-dimensional lens space, $W_{\xi}$ is a simply connected symplectic 
manifold, with first Chern number $N_{W_{\xi}}=p$. It is monotone
for $\xi < 0$ and the lens space $L^{2n+1}_p$ can be embedded in
$W_{\xi}$ as a Lagrangian submanifold.\\ 

\textbf {Example 3:} For $n \geq 3$, $\widetilde{G}_2(\rr^{n+2})$ is
simply connected and the Euler class of the bundle $V_2(\rr^{n+2}) \rightarrow
\widetilde{G}_2(\rr^{n+2})$ is equal to $-1$
times the generator. The symplectic cut $W_{\xi}$ is hence a simply
connected manifold with first Chern number $N_{W_{\xi}}=1$. For $\xi <
0$, $W_{\xi}$ is a monotone symplectic manifold, in which the Stiefel
manifold $V_2(\rr^{n+2})$ is embedded as a Lagrangian submanifold.

For the low dimensions $n=1$ and $n=2$, the Grassmannian 
$\widetilde{G}_2(\rr^{n+2})$ and the symplectic cut are simply
connected with Chern numbers: for $n=1$,
$N_{W_{\xi}}=2$ and for $n=2$, $N_{W_{\xi}}=1$.\\

\section{Monotonicity of the Lagrangian submanifold $V$ in $W$}
\label{sec:MonotonieV}

We still assume that $V \rightarrow B$ is a principal circle bundle
and that $W_{\xi}$ is the symplectic cut of the cotangent bundle of
$V$ at the level $\xi$ for a negative real number $\xi$ (in order to
simplify the notations, we will omit all the indices $\xi$ in the following).\\

Recall (see \cite{MR1223659}) that a Lagrangian submanifold $L$ of a 
symplectic manifold~$(W,\omega)$ is monotone if there exist a positive 
constant $K_L$ such that for all $v$ in $\pi_{2}(W,L)$, 
$$\int v^\star \omega = K_L \, \mu_{L} (v)$$
where $\mu_{L} : \pi_{2}(W,L) \rightarrow \zz$ denotes the Maslov class of 
$L$ in $W$ (we will recall the definition of the (relative) Maslov class in the
proof of Lemma \ref{lem:Maslov-gen}).

Moreover, if $L$ is a Lagrangian submanifold of a symplectic manifold $W$, it 
is known (see \cite{MR1223659}) that $W$ must be monotone and that the
constants $K_W$ and $K_L$ satisfy:
$$K_W = 2 K_L.$$
We will also denote by $N_L$ the Maslov number of $L$, namely the
nonnegative generator of the subgroup $\langle \mu_{L} , \pi_2(W,L)
\rangle$ of $\zz$.

\begin{rem} \label{rem:NLetNW}
  \begin{description}
  \item[(i)] If a disc $v$ of $\pi_{2}(W,L)$ is the image of a sphere
    $w$ of $\pi_{2}(W)$, we have the equalities: $\mu_{L}(v)=2 c_{1}(w)$
    and $\int v^\star \omega = \int w^\star \omega$ (see for instance \cite{MR1223659}). 
    In particular, $N_{L}$ always divides $2 N_{W}$,
  \item[(ii)] In the case of a simply connected Lagrangian submanifold
    $L$, the map $\pi_{2}(W) \rightarrow \pi_{2}(W,L)$ is
    surjective. Consequently, if $W$ is monotone, $L$ is monotone and 
    $N_{L} = 2 N_{W}$.
  \end{description}
\end{rem}

In \cite{Birannew}, Biran used a remark similar to (ii) for
submanifolds $L$ such that $H_{1}(L;\zz)$ is of $q$-torsion (i.e. for
each $\alpha \in H_{1}(L;\zz), q \alpha = 0$). Here, we will use the
following: 

\begin{lem}
\label{lem:Biran}
If $L$ is a Lagrangian submanifold of a monotone symplectic manifold
$(W,\omega)$ and if $\pi_{1}(L)$ is of $q$-torsion (namely $a^{q} = 1$
for any $a$ of $\pi_{1}(L)$ with $q \neq 0$), then $L$ is monotone in
$W$ and $2N_{W}$ divides $q N_{L}$.
\end{lem}

\begin{proof}
  As $W$ is monotone, there exists a positive constant
  $K_W$ such that for any $w$ in $\pi_{2}(W)$,
\begin{equation}
  \int w^\star \omega  = K_L \, \mu_{L} (v).
\label{eq:B-rel-mon}
\end{equation}

Let $v$ be an element of $\pi_{2}(W,L)$. The boundary of $v$ is an
element of $\pi_{1}(L)$ and thus if $L$ is of $q$-torsion, 
$(\partial v)^q = \partial (v^q) = 1$ in $\pi_{1}(L)$. This means that $v^q$
is a sphere $w$ of $W$.

We have then the following relations:
\begin{equation}
  2 c_{1}(w) = \mu_{L} (v^q) = q \mu_{L} (v)
\label{eq:B-c1}
\end{equation}
and
\begin{equation}
   \int w^\star \omega  =  \int (v^q)^\star \omega = q \int v^\star \omega .
\label{eq:B-omega}
\end{equation}
The relation (\ref{eq:B-rel-mon}) implies:
$$q \int v^\star \omega = 2 K_W q \mu_{L} (v)$$
and hence the monotonicity of $L$ in $W$, since $q$ is nonzero.

Moreover by relation (\ref{eq:B-c1}), for all $v$ in $\pi_{2}(W,L)$, $2 N_{W}$
divides $q \mu_{L} (v)$ and hence $2 N_{W}$ divides $q N_{L}$.
\end{proof}

Remark \ref{rem:NLetNW} and Lemma \ref{lem:Biran} are useful to study
our examples:\\

\textbf {Example 0:} For the action of $\sph^{1}$ on itself by
multiplication, we know that 
$$\pi_{2}(W,\sph^{1}) = \pi_{2}(\cc,\sph^{1}) \simeq \pi_{1}(\sph^{1}) \simeq
\zz.$$ 
Moreover, the Maslov class of the disc centered at the origin and of 
radius $\sqrt{-2 \xi}$ (a generator of $\pi_{2}(W,\sph^{1})$) is $2$ 
(it is a disc in $\cc$) and its area is $- 2 \pi \xi$. This 
means that the area of the disc is equal to $- \pi \xi$ times its 
Maslov class and that $\sph^{1}$ is monotone in $\cc$.\\

\textbf {Examples 1 and 3:} In the case of the Hopf bundle (or of the 
Stiefel manifolds), as the sphere $V=\sph^{2n+1}$ for $n \geq 1$ (or the 
Stiefel manifold $V = V_2(\rr^{n+2})$ for $n \geq 3$) is simply 
connected, $V$ is monotone in $W$ and $N_{V}=2N_{W}=~2$. These
manifolds $V$ are even $2$-connected, hence the long exact sequence
of the pair $(W,V)$ gives an isomorphism between $\pi_{2}(W,V)$ and
$\pi_{2}(\pnc) \simeq \zz$ (or 
$\pi_{2}(\widetilde{G}_2(\rr^{n+2})) \simeq \zz$).\\

\textbf {Example 2:} If $V$ is the lens space $L_{p}^{2n+1}$, its
fundamental group is $\zz/p$ and in particular, $\pi_{1}(V)$ is of
$p$-torsion. By Lemma~\ref{lem:Biran}, $V$ is monotone
in $W$ and $2 N_{W}$ divides $p N_{V}$, namely $2 p$ divides $p N_{L}$
and hence $2$ divides $N_{V}$.\\

Let us prove that Example 0 essentially describes what happens in
general for the Lagrangian submanifold $V$ in the symplectic cut of
its cotangent bundle, namely:

\begin{thm}
\label{thm:Vmonotone}
Let $\pi : V \rightarrow B$ be a principal circle bundle and let $W$ 
be the symplectic cut of $\cotg$ at the level $\xi$ with $\xi < 0$.\\
Then, $\pi_{2}(W,V)$ is isomorphic to $\zz$, the Lagrangian
submanifold $V$ is monotone in $W$ with monotonicity constant $-\pi \xi$ 
and its Maslov number is $2$.
\end{thm}

We begin with the computation of $\pi_{2}(W,V)$:

\begin{lem} 
\label{lem:pi2mon}
Let $\pi : V \rightarrow B$ be a principal circle bundle and let $W$ be 
the symplectic cut of $\cotg$ at the level $\xi$, with $\xi < 0$.\\
Then $\pi_{2}(W,V)$ is isomorphic to $\zz$, a generator $\phi$ being 
given by the image of the disc of $\cc$ centered at the origin and of
radius $\sqrt{-2 \xi}$ in a fibre of the complex line bundle 
$W \rightarrow Q$. 
\end{lem}

\begin{proof}
In order to prove that $\pi_{2}(W,V)$ is isomorphic to
$\pi_{1}(\sph^{1})$, we connect the long exact sequence of the pair
$(W,V)$ to the long exact sequence associated with the circle bundle $V
\rightarrow B$:
$$
\xymatrix{ \pi_{2}(V) \ar[r] \ar[d]^-{=} & \pi_{2}(B) \ar[r]
  \ar[d]^-{\simeq} &
  \pi_{1}(\sph^{1}) \ar[r]  \ar@{.>}[d]^-{\Phi} & \pi_{1}(V) \ar[r] \ar[d]^-{=} &  \pi_{1}(B) \ar[r] \ar[d]^-{\simeq}  &  0\\
  \pi_{2}(V) \ar[r] & \pi_{2}(W) \ar[r] & \pi_{2}(W,V) \ar[r] &
  \pi_{1}(V) \ar[r] & \pi_{1}(W) \ar[r] & 0 }
$$
and we define a map $\Phi: \pi_{1}(\sph^{1}) \rightarrow \pi_{2}(W,V)$
that makes the above diagram commute. The map $\Phi$ is given by
associating to a generator of $\pi_{1}(\sph^{1})$ a disk $\phi$, image
of the disc of $\cc$ centered at the origin and of radius $\sqrt{-2
  \xi}$ in a fibre of the complex line bundle $W \rightarrow Q$. 
(see Appendix \ref{sec:Append-Constr-map-Phi} for more details).

By the Five-lemma, we conclude that $\Phi$ is an isomorphism (hence, 
$\pi_{2}(W,V)$ is isomorphic to $\zz$) and that the image of 
$\pi_{2}(W) \rightarrow \pi_{2}(W,V)$ is equal to the image of 
$\pi_{2}(B) \rightarrow \pi_{1}(\sph^{1})$, namely 
$N_{e} \zz = \langle e , \pi_{2}(B) \rangle$. \end{proof}

We can determine the Maslov number of $V$ in $W$ thanks to the
following commutative diagram:
$$
\xymatrix{
\pi_{2}(W) \ar[r]^-{\Phi \circ e} \ar[d]^-{=} &  \pi_{2}(W,V) \ar[r]
\ar[d]^-{\mu_{V}} & \pi_{1}(V)  \ar[d]^-{\gamma_{V}} \ar[r] &
\pi_{1}(W) \ar[d]^-{=} \\
\pi_{2}(W) \ar[r]^-{2 c_{1}(W)} &  \pi_{1}(\Lambda_{n}) \ar[r] & \pi_{1}(\Lambda(W)) \ar[r] & \pi_{1}(W)
         }
$$
where $\Lambda_n$ (for $n$ half the dimension of $W$) is the Grassmannian 
of the Lagrangian subspaces of $\rr^{2n}$ and $\Lambda(W)$ is the fibre 
bundle associated to the vector bundle $TW \rightarrow W$ with typical fibre $\Lambda_n$.

We have $\mu \circ (\Phi \circ e) = 2 c_{1}(W)$ and hence if the Euler class $e$ 
does not vanish on $\pi_{2}(B)$ (in particular, $\pi_{2}(B)$ cannot 
be trivial), then 
$\mu_{V} : \pi_{2}(W,V) \simeq \zz \rightarrow \pi_{1}(\Lambda_{n}) \simeq \zz$
 is the multiplication by $2$ or $-2$ and the Maslov number of $V$ is 
$N_{V} = 2$.\\

In particular, if $e$ is not trivial on  $\pi_{2}(B)$, we can say that
the Maslov class of the generator of $\pi_{2}(W,V)$ is $2$ or
$-2$. Unfortunately, if we want to conclude for the monotonicity of
$V$, we must find out the sign of this class. We use one more time
Example 0.

\begin{lem}
\label{lem:Maslov-gen}
  The Maslov class in $W$ of the generator $\phi$ is $2$.
\end{lem}

\begin{proof} 
By definition, we get the Maslov class $\mu_{V}^{W}$ of a disc 
$v:(D^{2},\sph^{1}) \rightarrow (W,V)$ by trivialising the tangent bundle
$TW$ on the disc:
$$
\xymatrix{
  D^{2} \times \rr^{2n} \simeq v^{\star} TW \ar[r] \ar[d] & TW \ar[d] \\
  D^{2} \ar[r]^-{v} & W }
$$
As the restriction $v_{|\sph^{1}}$
takes its values in $V$, $(v_{|\sph^{1}})^{\star}TV$ defines a loop in
$TV$ and after trivialisation of $v^{\star}TW$, a loop in the
Grassmannian $\Lambda_{n}$ of the linear Lagrangian subspaces of $\rr^{2n}$. 
The class $\mu_{V}(v)$ is then the (ordinary) Maslov
class of this loop of $\Lambda_{n}$.

In our case, the disc $\phi$ lies in the fibre of $W \rightarrow
Q$. We describe $TW$ in the neighbourhood of this fibre. Let
$\mathcal{U}$ be a contractile neighbourhood of $b_{0}$ in $B$ so that
the bundle $V \rightarrow B$ can be trivialised on $\mathcal{U}$:
$$
\xymatrix{
  \! V_{|\mathcal{U}} \ar[rr]^-{~} \ar[rd]^-{\pi} & & \mathcal{U} \times \sph^{1} \ar[ld] \\
  & \mathcal{U} & }
$$
On the open subset $V_{|\mathcal{U}}$ of $V$, the restriction of $\pi$
is a trivial bundle. We are therefore in the case of Example 0. This
means that $W_{|\mathcal{U}}$ ($= (p \circ q)^{-1}(\mathcal{U})$) is
symplectomorphic to $T^{\ast}\mathcal{U} \times \cc$ and $V_{|\mathcal{U}}$ is embedded 
in $T^{\ast}\mathcal{U} \times \cc$ as the product of the zero section 
of $T^{\ast}\mathcal{U}$ and of the circle centered at the origin and
of radius $\sqrt{-2 \xi}$. The Maslov class of the disc $\phi$ is
hence equal to the Maslov class of the loop $u \mapsto \sqrt{-2 \xi}
u$ in $\cc$, namely $2$. \end{proof}

\textit{Proof of Theorem \ref{thm:Vmonotone}.} 
Thanks to the symplectomorphism between $\cc$ and the fibre (see
Appendix \ref{sec:Append-Constr-map-Phi}), we also
know the area of the disc $\phi$: it is equal to the area of the disc 
of $\cc$ with center at the origin and radius $\sqrt{-2 \xi}$, 
namely~$-2 \pi \xi$. As in the case of $T^{\ast}{\sph^{1}}$, 
$$\int v^{\star} \omega = - \pi \xi \; \mu_{V}^{W}(v)$$
if $v$ is the disc $\phi$ and thus for all $v$ of 
$\pi_{2}(W,V)$.\hfill \qedsymbol\\

In the next section, we are interested in the other monotone 
Lagrangian embeddings in a symplectic cut.

\section{Monotone Lagrangian submanifolds in symplectic cuts}
\label{sec:GeneralLag}

As annonced in Introduction, we will apply periodicity 
results of Seidel~\cite{MR1765826} to Lagrangian submanifolds 
in the symplectic cut. We recall here some 
general results on this periodicity.

\subsection{Seidel's periodicity theorem}
\label{sec:Seidel}

To begin with, let us state the theorem of Seidel~\cite{MR1765826} (see 
also~\cite{MR2191628}):

\begin{thm}[Seidel] 
\label{thm:S-per}
Let $(W,\omega)$ be a monotone geometrically bounded symplectic 
manifold endowed with a Hamiltonian circle action.
Let $w$ be the sum of the weights of the linearised action at a fixed 
point.\\
Assume that there exists a positive integer $N$ such that 
$2c_{1}(W,\omega)$ is mapped on zero in $H^{2}(W;\zz/N)$.\\
Let $L$ be a monotone Lagrangian submanifold of $W$ (with $N_{L} \geq
2$) such that the modulo $N$ Maslov class 
$\left({ \gamma_{L}^{W} }\right)^{\star}(\mu^{N})$ of $L$
vanishes in $H^{1}(L;\zz/N)$.\\
Then the Floer cohomology of $L$ is (absolutely) graded by $\zz/N$ and
is periodic of period $2 w$.
\end{thm}

Let us recall how the Gauss map $\gamma_{L}^{W}$ and the reduced 
Maslov class $\mu^{N}$ are defined .

If $(W,\omega)$ is a symplectic manifold, we consider as 
in~Section~\ref{sec:MonotonieV} the bundle $\Lambda(W) \rightarrow W$ 
associated with the tangent bundle whose fibre at a point is the 
set of linear Lagrangian subspaces of the tangent space to $W$ at 
that point.
If $L$ is a Lagrangian submanifold of $W$, we denote by 
$\gamma_{L}^{W}: L \rightarrow \Lambda(W)$ the Gauss 
map that associates to a point of $L$ the tangent space of $L$ 
at that point. 

If $(W,\omega)$ is a symplectic manifold such that 
$2 c_{1}(W,\omega) = 0$ in $H^{2}(W;\zz/N)$ for some 
integer $N$, there exists (see for example \cite{MR1765826}) a
Maslov class $\mu^{N}$ in $H^{1}(\Lambda(W);\zz/N)$ that extends
the modulo $N$ reduction of the ordinary Maslov class on each fibre of
the bundle $\Lambda(W) \rightarrow W$. Note that when $W$ is a cotangent 
bundle, one can take $N = 0$ and 
$\left({ \gamma_{L}^{W} }\right)^{\star}(\mu^{N})$ is then the usual 
Maslov class of a Lagrangian submanifold $L$ in the cotangent bundle.\\

\begin{rem}
\label{rem:simplificationcasSC}
  If $W$ is a simply connected symplectic manifold, the assumption  
``$2c_{1}(W,\omega)$ is mapped on zero in $H^{2}(W;\zz/N)$'' can be
replaced by ``$N$ divides $2 N_{W}$'' and the assumption 
``$\left({ \gamma_{L}^{W} }\right)^{\star}(\mu^{N})$ vanishes in 
$H^{1}(L;\zz/N)$'' by ``$N$ divides $N_{L}$''.
\end{rem}

We recall now how the weights of a Hamiltonian circle action are defined.

If $(W, \omega)$ is endowed with a Hamiltonian circle action, this
action defines a map of $\sph^{1}$ in 
the group of symplectic diffeomorphisms of $W$: 
$$
\begin{array}{cccc}
\rho: & \sph^1 & \longrightarrow & \Sp(W)\\
      &  u  &  \longmapsto    & (x \mapsto u \cdot x)
\end{array}
$$
At a fixed point $x$, the linearised actions $T_{x}\rho(u)$ (for 
$u \in \sph^{1}$) are endomorphisms of $T_{x}W$ whose matrices in a
common diagonalisation basis have the following form:
$$
\left( \begin{array}{ccc} u^{m_1} &        & 0  \\ 
                                  & \ddots &    \\
                             0    &        & u^{m_n} \end{array}
                            \right)
$$
The $m_{i}$'s are called the weights of the action.\\

\textbf {Example:} In \cite{MR1765826}, Seidel applied Theorem
\ref{thm:S-per} to the projective space $\pnc$. 
Together with the theorem of Oh, 
\begin{thm}[Oh, \cite{MR1389956}]
\label{thm:Oh}
Let $L$ be a compact monotone Lagrangian submanifold of a symplectic
manifold $(W,\omega)$ such that $N_{L} \geq 2$.\\
Denote by $HF(L,L)$ its Floer cohomology.\\
\begin{description}
\item[(i)] If $N_{L} \geq n+2$, then 
$\displaystyle HF(L,L) \cong \bigoplus_{k} H^{k}(L;\zz/2)$;
\item[(ii)] If $N_{L} = n+1$, then 
$\displaystyle HF(L,L) \cong \bigoplus_{k}   H^{k}(L;\zz/2)$ or 
$\displaystyle \bigoplus_{k \neq 0, n}   H^{k}(L;\zz/2)$.
\end{description}
\end{thm}

\noindent he proved that there is no simply connected Lagrangian
submanifold in $\pnc$ (for the precise statement, see 
\cite[Theorem 3.1]{MR1765826}).
Seidel used the following circle action on the projective space: 
$$
\begin{array}{ccc}
\sph^1 \times \pnc & \rightarrow & \pnc \\
  (u,[z])          & \mapsto     & [u z_0, z_1, ..., z_n]
\end{array}
$$
of Hamiltonian  
$$H([z]) = \frac{|z_0|^2}{\sum_{i=0}^{n} |z_i|^2}.$$
The fixed points are the points of the hyperplane $z_0 = 0$, for which
the sum of weights is $1$ and the point of homogeneous coordinates
$[1,0,...,0]$ for which the sum of weights is $-n$.\\

We can notice that, in this example, $N_W = n+1$ and that the classes of
the sums of weights are equal in $\zz/N_W$. Thus, we get for all the
fixed points the same period $2$ when the Floer cohomology is graded
by $\zz/N$ with $N$ dividing $N_W$. This is a general fact:

\begin{prop}
\label{prop:poids}
The class in $\zz/N_W$ of the sum of weights at a fixed point does not 
depend on the fixed point.
\end{prop}

This may be well known to specialists. We nevertheless include a
proof in Appendix \ref{sec:Append-Proof-Prop-poids}.\\

In the case of the symplectic cut of the cotangent bundle, we have:

\begin{prop} \label{prop:per}
The symplectic cut of the cotangent bundle $W_{\xi}$ is endowed with a
Hamiltonian circle action and the sum of the weights of the linearised 
action in a fixed point is equal to $1$. 
\end{prop}

\begin{proof}
The symplectic cut $W_{\xi}$ of the cotangent bundle is endowed with
the Hamiltonian circle action: 
  $$
\begin{array}{ccc}
\sph^1 \times W_{\xi} & \longrightarrow & W_{\xi} \\
  (u,[x,\varphi,z])&  \longmapsto    & [x, \varphi, uz]
\end{array}
$$
of Hamiltonian
$$
\begin{array}{cccc}
h: &    W_{\xi}   & \longrightarrow & \rr \\
   & [x,\varphi,z]& \longmapsto     & \frac{1}{2} |z|^2.
\end{array}
$$
As the action of $\sph^1$ is free on $\cotg$, the fixed points of this
action are the points such that $z=0$, namely the points lying in
$Q_{\xi}$.

At a point $(x,\varphi)$ of $Q_{\xi}$, the action of $\sph^1$ can be
linearised under the decomposition
$$T_{[x,\varphi,0]} W_{\xi} \simeq T_{[x,\varphi]} Q_{\xi} \oplus
(H^{-1}(\xi) \times_{\sph^1} \cc)_{[x,\varphi]}$$
by
$$
\left(
\begin{matrix}
  \mbox{id} & 0 \\
  0 & u
\end{matrix} \right)
$$
It implies that the sum of weights is $1$.
\end{proof}

\begin{rem}
This is (modulo $N_{W_{\xi}}$) the smallest nonzero sum of weights 
we can expect. When the assumptions of Theorem \ref{thm:S-per} are 
fulfilled, the Floer cohomology of a Lagrangian submanifold is 
$2$-periodic.
\end{rem}

We will use Seidel's theorem in the next section in order to get obstructions 
to the existence of Lagrangian submanifolds into a symplectic cut. 
To apply Seidel's theorem to the symplectic cut, we have to check the 
following : 

\begin{prop}
   The symplectic cut is geometrically bounded in the sense 
   of~\cite{MR1274934}.
\end{prop}

\begin{proof}
Recall (see~\cite{MR1274934}) that a symplectic manifold without 
boundary $(W,\omega)$ is geometrically bounded if there exists 
on~$W$ an almost complex structure $J$ and a complete Riemannian 
metric $g$ such that:
\begin{description}
\item[a)] $J$ is uniformly tamed by $\omega$, that is there exist 
strictly positive constants $\alpha$ and $\beta$ such that:
\begin{equation*}
   \label{eq:geobornalpha}
   \omega(X,JX) \geq \alpha \; g(X,X)
\end{equation*}
and
\begin{equation}
   \label{eq:geobornbeta}
   |\omega(X,Y)| \leq \beta \; \|X\|_g \|Y\|_g
\end{equation}
for all $X,Y \in TW$.

\item[b)] There exist an upper bound for the sectional curvature 
of $(W,g)$ and a stricly positive lower bound for the injectivity 
radius of $(W,g).$
\end{description}

We also recall that the cotangent bundle $(\ctg V, \omega_{\can})$ of 
a closed manifold $V$ is geometrically bounded. To see it, we can 
choose a metric $g$ on $\ctg V$ induced by a Riemannian metric 
on $V$ and an $\omega_{\can}$-tame almost complex structure $J$ 
homogeneous with respect to uniform dilatations in the fibres. 
We can choose for example the almost complex structure induced by 
the Levi-Civit\`a connection (as in \cite{MR1001276} and \cite{MR1389956}).

The product of $(\ctg V , \omega_{\can})$ 
with $(\cc,\omega_{\std}, J_{\std}, g_{\std})$ is then also geometrically 
bounded.

When $V$ is endowed with a free circle action, one can choose a 
Riemannian metric on $V$ which is invariant by the action of the 
circle. The metric on the product $\ctg V \times \cc$ is then also 
invariant for the induced Hamiltonian action so that we get by 
restriction the upper bounds and lower bound of a) and b) on the 
level $\tH^{-1}(\xi)$. Taking the quotient, the metric and the almost 
complex structure induced on $W_{\xi}$ satisfy a) and b), so that $W_\xi$ is 
geometrically bounded.
\end{proof}

\subsection{Simply-connected embeddings}

For a monotone Lagrangian embedding of a simply connected manifold 
in a symplectic cut, we have the following topological obstructions:

\begin{thm}
\label{thm:nsc-spcut}
Let $B$ be a compact manifold of dimension $d-1$ with $d \geq 2$.\\
Let $e$ be an element of $H^{2}(B;\zz)$ such that $e$ is
not trivial on $\pi_{2}(B)$ and $2 e = 0$ in $H^{2}(B;\zz/N)$ for an
integer $N > 2$.\\ 
Let $V^{d} \rightarrow B^{d-1}$ be the principal circle bundle with
Euler class $e$.\\
We denote by $W$ a monotone symplectic cut of $\cotg$ (that is a symplectic 
cut at a negative level).\\
If $N > d+2$ or $N = d+2$ with $d$ odd, there is no compact and simply 
connected Lagrangian submanifold in $W$.\\
If $N = d+2$ with $d$ even, the $\zz/2$-cohomology groups of a compact
and simply connected Lagrangian submanifold in $W$ are isomorphic
to those of $\mathbb{C} \mathbb{P}^{d/2}$.
\end{thm}

\begin{proof}
  If $2 e = 0$ in $H^{2}(B,\zz/N)$ for some integer $N > 2$, $N$ 
  divides $2 N_{e}$, so that $2 N_{e} = 2 N_W > 2$.

  Let $L$ be a compact and simply connected Lagrangian submanifold 
  in~$W$.   By Remark~\ref{rem:NLetNW} (ii), its Maslov number in 
  the symplectic cut $W$ is $N_{L}= 2 N_{W}$, 
  and in particular $N_L \geq 2$. 
  Applying Theorem~\ref{thm:S-per}, the Floer homology of $L$ is well 
  defined, graded by $\zz/N$ and $2$-periodic.

  If moreover $N \geq d+2$, then $N_{L} \geq d+2$ and we can apply 
  Oh's theorem: the Floer homology of $L$ is isomorphic to the ordinary 
  cohomology of $L$ with   coefficient in $\zz/2$:
  $$HF(L,L) \cong \bigoplus_{i} H^i(L;\mathbb{Z}/2\mathbb{Z}).$$
  This leads to a contradiction with the $2$-periodicity if $N \geq d+3$ 
  (see Figure~\ref{fig:cerchomol}); and if $N = d+2$, the periodicity 
  implies   that the cohomology groups of $L$ are isomorphic to the 
  cohomology groups of $\mathbb{C} \mathbb{P}^{d/2}$ for even 
  integer $d$ and are trivial if $d$ is odd, which gives an other 
  contradiction.

  \begin{figure}[htbp]
  \begin{center}
   \psfrag{2m = 0}{$N \equiv 0$}
   \psfrag{Z/2}{$\zz/2$}
   \psfrag{0}{$0$}
   \psfrag{d}{$d$}
   \includegraphics{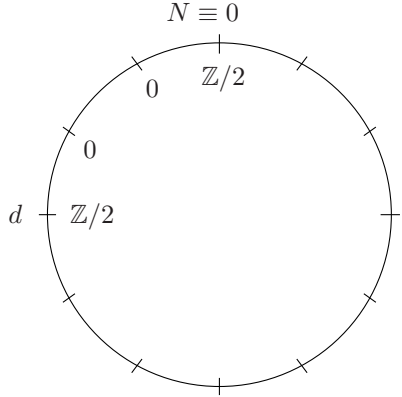}
    \caption{the Floer homology of $L$ when $N_{L} \geq d+2$}
    \label{fig:cerchomol}
  \end{center}
\end{figure}
\end{proof}

\subsection{Lagrangian submanifolds of the cotangent bundle}
\label{sec:Ss-var-ds-cotg}

Let $L$ be a compact Lagrangian submanifold in~$\cotg$. Then, for 
$|\xi|$ large enough, it lies in $H^{-1}(\xi;+\infty)$ and hence 
can be embedded as a Lagrangian submanifold of the symplectic cut 
$W_{\xi}$. Such a $\xi$ being now fixed, we will omit the 
indices~$\xi$ in the following.

\subsubsection{Monotonicity of a Lagrangian submanifold of the cotangent 
bundle into the symplectic cut}
\label{sec:mon-lag-ctg-into-spcut}

In order to understand under which conditions $L$ can be monotone in
$W$, we need to know its Maslov class in $W$. In the following, 
we will relate this class to the Maslov class of $L$ in~$\cotg$.

Remember that the cotangent bundle $\cotg$ of $V$ is a symplectic manifold the
first Chern class of which satisfies $2 c_{1}(\cotg) = 0$ in
$H^{2}(\cotg;\zz)$. As already noticed in Section~\ref{sec:Seidel}, 
this implies that there exists a Maslov class $\mu^{\cotg} \in H^{1}(\Lambda(\cotg);\zz)$ 
that extends the ordinary Maslov class on
each fibre of $\Lambda(\cotg) \rightarrow \cotg$. The Maslov class 
$\mu_{L}^{\cotg} \in H^{1}(L;\zz)$ of a Lagrangian submanifold 
$L$ in $\cotg$ can then be defined as the pullback: 
$$\mu_{L}^{\cotg} = (\gamma_{L}^{\cotg})^{\star}\mu^{\cotg}.$$

The formula connecting the Maslov class of $L$ in $W$ to that of $V$
in $W$ and that of $L$ in $\cotg$ is the analogue of the formula given by 
Lalonde and Sikorav in \cite{MR1090163}. We use their notations and
denote by $i$ the embedding of $L$ in $\cotg$ and $f$ the composition
of this embedding and of the projection of $\cotg$ on $V$.
$$
\xymatrix{
  L \ar[r]^-{i} \ar[rd]^-{f} & \cotg \ar[d]^-{p_{V}} \\
                             & V }
$$ 
We use the long homotopy exact sequence associated with the triple 
$(W,\cotg,L)$:
$$\pi_{2}(\cotg,L) \longrightarrow \pi_{2}(W,L) \longrightarrow
\pi_{2}(W,\cotg) \longrightarrow \pi_{1}(\cotg,L).$$
The projection is a homotopy equivalence between $\cotg$ and $V$, so
that $\pi_{2}(W,\cotg)$ is isomorphic to $\pi_{2}(W,V)$ and the 
composition of this isomorphism with the map 
$\pi_{2}(W,L) \rightarrow \pi_{2}(W,\cotg)$ 
is a map that will be denoted $f_{\star}$ (it is actually 
$(\ident,f)_{\star}$).

We can now write the formula:

\begin{lem} If $\partial:\pi_{2}(W,L)\rightarrow \pi_{1}(L)$ denotes the
boundary map,
\begin{equation}
  \label{eq:formuleLS}
  \mu_{L}^{W} = \mu_{V}^{W} \circ f_{\star} + \mu_{L}^{\cotg} \circ
\partial.
\end{equation}
\end{lem}

\begin{proof}
Let $w: (D^{2},\sph^{1}) \rightarrow (W,L)$ be a disc with boundary on $L$. 
The map $\pi_{2}(W,L) \rightarrow \pi_{2}(W,\cotg)$ associates to the 
class of $w$ seen as a disc with boundary in $L$ the class of $w$ 
seen as a disc with boundary in $\cotg$. Thanks to the retraction in 
the fibres of the projection $p_{V}: \cotg \rightarrow V$, we have a  
homotopy $h$ between $\partial w$ and $p_{V}(\partial w)$ and we can 
describe the image of the class of $w$ by the isomorphism 
$\pi_{2}(W,\cotg) \rightarrow \pi_{2}(W,V)$ as the class of the disc 
$\widetilde{w}= w \# h$ obtained by gluing $w$ and $h$ (see Figure 
\ref{fig:ssvarlagcotg}), so that $f_{\star}[w]=[\widetilde{w}]$.

\begin{figure}[htbp]
  \begin{center}
   \psfrag{V}{$V$}
   \psfrag{T*V}{$\cotg$}
   \psfrag{L}{$L$}
   \psfrag{Q}{$Q$}
   \psfrag{w}{$w$}
   \psfrag{h}{$h$}
   \psfrag{f}{$f$}
   \includegraphics{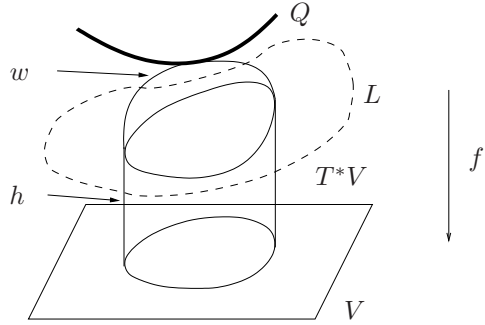}
    \caption{$\widetilde{w}$ obtained by gluing $w$ and $h$}
    \label{fig:ssvarlagcotg}
  \end{center}
\end{figure}

Since $h$ takes value in $\cotg$, the difference  
$$ \mu_{L}^{W}(w) - \mu_{V}^{W} \circ f_{\star}(w) = \mu_{L}^{W}(w) -
\mu_{V}^{W}(\widetilde{w})$$
is simply the Maslov class in $\cotg$ of the boundary of $h$ (seen as
a map of $[0,1] \times [0,1]$ in $\cotg$) namely the difference
between the Maslov classes in $\cotg$ of the boundaries of $w$ and of 
$\widetilde{w}$. And $\widetilde{w}$ is a loop in the zero section $V$ 
of $\cotg$, therefore its Maslov class vanishes and we have:
$$\mu_{L}^{W}(w) - \mu_{V}^{W}(\widetilde{w}) = \mu^{\cotg}(\partial w).$$
\end{proof}

In order to prove the monotonicity of $L$, we need a similar formula
for the class of the symplectic form. With the notations of the
previous proof, this relation can be written:
\begin{equation}
  \label{eq:LSomega}
  \int_{D^{2}} w^{\star} \omega_{W}= \int_{D^{2}}
  \widetilde{w}^{\star}\omega_{W} + \int_{\sph^{1}} (\partial
  w)^{\star}\lambda 
\end{equation}
where $\lambda$ is the Liouville $1$-form of $\cotg$.\\

Indeed, we have $\widetilde{w}=w \# h$, hence 
\begin{eqnarray*}
\int_{D^{2}} \widetilde{w}^{\star} \omega_{W}
& = & \int_{D^{2}} w^{\star}\omega_{W} + \int_{[0,1]\times[0,1]}
h^{\star}\omega_{W}\\
& = & \int_{D^{2}} w^{\star}\omega_{W} + \int_{\sph^{1}} (p_{V}(\partial w))^{\star}\lambda - \int_{\sph^{1}} (\partial w)^{\star}\lambda 
\end{eqnarray*}
But $p_{V}(\partial w)$ is a loop in $V$, and $V$ is exact in $\cotg$, 
hence the integral $\int_{\sph^{1}} (p_{V}(\partial w))^{\star}\lambda$ 
is zero.\\

\begin{cor}
\label{cor:exactedsctgdsspcut}
Any exact compact Lagrangian submanifold with zero Maslov class in the
cotangent bundle of $V$ can be embedded as a monotone Lagrangian
submanifold in a symplectic cut of $\cotg$ (for some negative level
$\xi$). 
\end{cor}

\begin{proof}
If we restrict to the case of a Lagrangian submanifold $L$ of the 
cotangent bundle that is exact and with zero Maslov class, the 
monotonicity of $V$ (Theorem~\ref{thm:Vmonotone}) implies the 
monotonicity of $L$ in $W$ as (\ref{eq:formuleLS}) and 
(\ref{eq:LSomega}) can then be written: for any disc 
$w:(D^{2},\sph^1) \rightarrow (W,L)$, 
\begin{eqnarray} \label{eq:LS-ex}
  \mu_{L}^{W}([w]) & = & \mu_{V}^{W}(f_{\star}[w]) = \mu_{V}^{W}([\widetilde{w}]), \\
  \int_{D^{2}} w^{\star} \omega_{W} & = & \int_{D^{2}}
  \widetilde{w}^{\star}\omega_{W}. \nonumber
\end{eqnarray}
\end{proof}

\subsubsection{Exact compact Lagrangian submanifolds of the cotangent bunble}

\subsubsection*{Verification of the assumptions of Seidel's periodicity theorem}

We verify now if Theorem \ref{thm:S-per} can be applied in the 
symplectic cut to the compact exact Lagrangian submanifolds of 
the cotangent bundle with vanishing Maslov class.

We have recalled in Section \ref{sec:mon-lag-ctg-into-spcut} why we can define 
a Maslov class $\mu^{\cotg} \in H^{1}(\Lambda(\cotg);\zz)$. We will
denote by 
$\left({\mu^{\cotg}}\right)^{N} \in H^{1}(\Lambda(\cotg);\zz/N)$ the 
modulo~$N$ reduction of $\mu^{\cotg}$.
As the embedding of $H^{-1}(\xi,+\infty)$ in $W_{\xi}$ is
symplectic, we have an embedding 
$$i_{\Lambda}: \Lambda(H^{-1}(\xi,+\infty)) \rightarrow \Lambda(W)$$
and if $2 c_{1}(W_{\xi})=0$ in $H^{2}(W_{\xi};\zz/N)$ (namely, if
$2e=0$ in $H^{2}(B;\zz/N)$), then  
$$\left({\mu^{\cotg}}\right)^{N}_{|\Lambda(H^{-1}(\xi,+\infty))} =
i_{\Lambda}^{\star}\mu^{N}.$$

Thus, if $L$ is a Lagrangian submanifold of $H^{-1}(\xi,+\infty)$,
\begin{eqnarray*}
  \left(\gamma_{L}^{W}\right)^{\star}\mu^{N} & = & \left (i_{\Lambda}
    \circ \gamma_{L}^{\cotg}\right)^{\star}\mu^{N} \\
 & = & \left(\gamma_{L}^{\cotg}\right)^{\star}
  \left ( i_{\Lambda} \right)^{\star}\mu^{N} \\
 & = & \left(\gamma_{L}^{\cotg}\right)^{\star}
 \left({\mu^{\cotg}}\right)^{N}_{|\Lambda(H^{-1}(\xi,+\infty))}
\end{eqnarray*}
Hence, $\left(\gamma_{L}^{W}\right)^{\star}\mu^{N}$ is the modulo $N$ 
reduction of the Maslov class 
$\left(\gamma_{L}^{\cotg}\right)^{\star} {\mu^{\cotg}}$ of $L$ in 
$\cotg$.

In particular, if $L$ is exact,
$\left(\gamma_{L}^{W}\right)^{\star}\mu^{N}$ is zero in 
$H^{1}(L;\zz/N)$.\\

\subsubsection*{Simply connected Lagrangian submanifolds}

Simply connected Lagrangian submanifolds are exact, in particular  
Theorem~\ref{thm:nsc-spcut} applies, so that we have 

\begin{thm}
\label{thm:nsc}
Let $B$ be a compact manifold of dimension $d-1$ with $d \geq 2$.\\
Let $e$ be an element of $H^{2}(B,\zz)$ such that $e$ is not trivial on 
$\pi_{2}(B)$ and $2 e = 0$ in $H^{2}(B;\zz/N)$ for some integer $N > 2$.\\ 
Let $V^{d} \rightarrow B^{d-1}$ be a principal circle bundle with Euler 
class $e$.\\
If $N > d+2$ or $N = d+2$ with $d$ odd, there is no compact simply 
connected Lagrangian submanifold in $\cotg$.\\
If $N = d+2$ with $d$ even, the cohomology groups with coefficients in 
$\zz/2$ of a compact and simply connected Lagrangian submanifold of $\cotg$ 
are isomorphic to those of $\mathbb{C} \mathbb{P}^{d/2}$.
\end{thm}

\subsubsection*{The case $B$ simply connected}

Lalonde et Sikorav proved in \cite{MR1090163} that if $L$ is an 
exact Lagrangian submanifold of $\cotg$, then the index of 
$f_{\star}(\pi_{1}(L))$ in $\pi_{1}(V)$ is finite. In the case $B$ 
(and hence $W$) is simply connected, Seidel's theorem 
(its simplified statement, see Remark~\ref{rem:simplificationcasSC}) 
gives:

\begin{thm}
\label{thm:cassc}
Let $V^{d} \rightarrow B^{d-1}$ be a principal circle 
bundle with base a compact and simply connected manifold $B$ and with
nonzero Euler class.\\
Let $L$ be an exact compact Lagrangian submanifold in $\cotg$. We assume 
that its Maslov class in the cotangent bundle is zero.\\
Then, if $m$ denotes the index of $f_{\star}(\pi_{1}(L))$ in
$\pi_{1}(V)$,
\begin{description}
\item [(i)] $m$ divides $N_{e}$ where $N_{e}$ the nonnegative generator of 
$\langle e , \pi_{2}(B) \rangle \subset \zz$;
\item [(ii)] $2m \leq d+2$.
\end{description}
\end{thm}

For the proof of Theorem~\ref{thm:cassc}, notice first that:

  \begin{rems}
  \label{rem:pi1V}
  \begin{description}
 \item 1) Let $V^{d} \rightarrow B^{d-1}$ be a circle bundle with base a 
    compact and simply connected manifold $B$.\\
    If $N_{e}$ is the nonnegative generator of the subgroup 
    $\langle e , \pi_{2}(B) \rangle$ of $\zz$, then $\pi_1(V)= \zz/N_e$. 
    In particular, $\pi_1(V)$ is cyclic.
 \item 2) If $L$ is a compact exact Lagrangian submanifold of $\cotg$ with 
    $B$ simply connected and if $m$ is the index of 
    $f_{\star}(\pi_{1}(L))$ in $\pi_{1}(V)$, 
    then the quotient
    $\pi_{1}(V)/f_{\star}(\pi_{1}(L)) \cong \zz/m$ is also cyclic.
\end{description}
  \end{rems}

\noindent
To see 1), one can compute $\pi_{1}(V)$ with the long exact sequence associated 
to the fibre bundle $V \rightarrow B$:
$$\pi_{2}(B) \longrightarrow \pi_{1}(\cerc) \longrightarrow \pi_{1}(V)
\longrightarrow \pi_{1}(B)=\{0\}.$$
As the image of the map $\pi_{2}(B) \rightarrow \pi_{1}(\cerc)$ 
is $N_{e} \zz$, $\pi_{1}(V) \cong \zz / N_{e}$ (so that we could have 
seen in this case that $V$ is monotone in $W$ thanks to 
Remark~\ref{lem:Biran}).\\

\noindent
\textit {Proof of Theorem~\ref{thm:cassc}:}
By Formula (\ref{eq:LS-ex}), $\mu_{L}^{W} = \mu_{V}^{W} \circ f_{\star}$.
As the image of $f_{\star}: \pi_2(W,L) \rightarrow \pi_2(W,\cotg)$ 
is $m \zz$ (because the index of $f_{\star}(\pi_{1}(L))$ 
in $\pi_{1}(V)$ is $m$) and $\mu_{V}^{W}$ is the multiplication by $2$, 
the Maslov number of $L$ is $N_{L}=2m$.
 
We also know that $N_{L}$ divides $2 N_{W} = 2 N_{e}$ 
(Remark~\ref{rem:NLetNW}), so $m$ necessarily divides $N_{e}$ 
and we have (ii).

Eventually, if $N_{L} \geq d+2$, then first, $N_{L} \geq 2$ so that we 
can consider the Floer homology of $L$ and secondly the assumptions of 
Oh's Theorem (Theorem~\ref{thm:Oh}) are fulfilled, so that the Floer 
homology is isomorphic to the cohomology of the submanifold:
$$HF(L,L) \cong \bigoplus_{i} H^i(L;\mathbb{Z}/2\mathbb{Z}).$$

By Theorem~\ref{thm:S-per} and Proposition~\ref{prop:per}, 
as $W$ is simply connected and as $2m$ divides $N_{L}$ and $2 N_{W}$, 
the Floer homology of $L$ is graded by $\zz/2m$:
$$HF^{i}(L,L)=H^{i}(L;\zz/2)$$
and by Proposition \ref{prop:per}, must be $2$-periodic. 
In the case when $N_{L} \geq d+3$, this periodicity is
in contradiction with Oh's theorem. 

\begin{cor}
\label{cor:Negrand}
  Under the assumptions of Theorem~\ref{thm:cassc}, if $2 N_{e} > d+2$, then
  $H^{1}(L;\zz/N_{e}) \neq \{0\}$. 
  In particular, $L$ cannot be simply connected (conclusion already contained 
  in Theorem~\ref{thm:nsc}).
\end{cor}

\begin{proof}
  If $2 N_{e} > d+2$, as $2m \leq d+2$ by Theorem~\ref{thm:cassc}, 
  the index $m$ of
  $f_{\star}(\pi_{1}(L))$  in $\pi_{1}(V)$ is strictly smaller than
  the cardinal $N_{e}$ of the group $\pi_{1}(V)$. Hence, the subgroup
  $f_{\star}(\pi_{1}(L))$ is not trivial, $f_{\star}$ is a nonzero
  homomorphism from $\pi_{1}(L)$ to $\zz/N_{e}$ and 
  $H^{1}(L;\zz/N_{e}) \neq \{0\}$. 
\end{proof}

Applied to the lens spaces, Theorem~\ref{thm:cassc} 
and Corollary~\ref{cor:Negrand}
give the following: 

\begin{prop}
If $L$ is a compact exact Lagrangian submanifold with zero Maslov
class in the cotangent bundle of a lens space $L_{p}^{2n+1}$, then the 
quotient $\pi_{1}(V)/f_{\star}(\pi_{1}(L))$ is cyclic and its 
cardinal, namely the index of $f_{\star}(\pi_{1}(L))$ in $\pi_{1}(V)$, 
is less than or equal to $n+1$ and divides $p$.\\
In particular, if $p$ is a prime number strictly greater than $n+1$,
then this index is equal to $1$.
\end{prop}

\begin{proof}
In the case of the lens spaces $V=L_{p}^{2n+1}$, $B=\pnc$ is simply
connected and $N_{e}=p$. The conclusion of Theorem~\ref{thm:cassc} 
can be written: $m$ divides $p$ and $2m \leq
2n+3$. 
Because of the parity, the case $2m = 2n+3$ cannot be realised, and
finally we have $m \leq n+1$.

In particular, if $p$ is a prime, $m$ is equal to $1$ or $p$ and if $p
> n+1$, then necessarily $m=1$.
\end{proof}

\begin{cor} 
Let $L$ be a Lagrangian submanifold in the cotangent bundle of a lens
space $L_{p}^{2n+1}$ for a integer $p$ such that $p > n+1$. Then 
$H^{1}(L;\zz/p) \neq \{0\}$.\\
In particular, $L$ cannot be simply connected nor with finite
fundamental group of cardinal $q$ with $\gcd(q,p)=1$
\end{cor}

But in the case the fundamental group of $V$ is cyclic, the following 
can be directly proved without using the symplectic cut of 
the cotangent bundle

  \begin{thm}
    \label{thm:referee}
    Let $V^{d} \rightarrow B^{d-1}$ be a circle bundle with basis a 
    compact simply connected manifold $B$.\\ 
    Let $L$ be a compact exact Lagrangian submanifold of $\cotg$.\\ 
    Then the map induced by $f$: 
    $$f_{\star}: \pi_{1}(L) \longrightarrow \pi_{1}(V)$$ 
    is surjective.
  \end{thm}

\begin{proof}
By Remark~\ref{rem:pi1V} 1), $\pi_1(V)= \zz/N_e$ is a quotient of 
$\pi_1(\sph^1)=\zz$.
  
To any abelian subgroup $m\zz/N_e \subset \zz/N_e$ ($m$ a divisor 
of $N_e$) is associated a finite cyclic covering with $m$ sheets 
$\widetilde{V} \rightarrow V$. The restriction of this covering 
to each fibre of the bundle $V \rightarrow B$ is just the 
path-connected covering with $m$ sheets 
of~$\sph^1$. In particular, if we lift the circle action $\rho(u)$ 
on $V$ (for $u \in \sph^1$) in a flow $\tilde{\rho}(u)$ on 
$\widetilde{V}$ (for $u \in \rr$), then $\tilde{\rho}(1)$, the 
flow at time $1$, is a generator of the covering group of 
$\widetilde{V} \rightarrow V$.

If we consider the associated covering 
$T^{\ast} \widetilde{V} \rightarrow \cotg$, the generator of its 
covering group (denoted by $\chi$) is isotopic to the identity 
through an Hamiltonian isotopy.
Now, assume that $j: L \rightarrow \cotg$ is an exact 
Lagrangian embedding such that the image $j_\star(\pi_1(L))$ is 
a subgroup $m \zz / N_e$ with $m>1$. This embedding can be 
lifted to  embedding 
$\tilde{j}: L \rightarrow T^\ast \widetilde{V}$, 
and it has the property that 
$\chi(\tilde{j}(L)) \cap \tilde{j}(L)$ is empty.
As we know that $\chi$ is Hamiltonian isotopic to the identity, 
we get a contradiction with the result of Gromov and Floer on the 
exact Lagrangian submanifolds. Consequently, $m=1$. 
\end{proof}

Let us consider now the monotone Lagrangian embedding of some particular 
manifold in the symplectic cut. We begin with the spheres.

\subsection{Monotone Lagrangian embeddings of spheres}

Theorem~\ref{thm:nsc-spcut} applies of course in the case of an 
embedding of sphere, but there exists even more obstructions in this case:

\begin{thm}
\label{thm:sphere}
Let $B$ be a compact manifold of dimension $d-1$
with $d \geq 2$.\\
Let $e$ be a nonzero element of $H^{2}(B;\zz)$ such that $e$ is
not trivial on $\pi_{2}(B)$ and $2 e = 0$ in $H^{2}(B;\zz/N)$ for an
integer $N > 2$.\\ 
Let $V^{d} \rightarrow B^{d-1}$ be a principal circle bundle with
Euler class $e$.\\
We denote by $W$ a monotone symplectic cut of $\cotg$. \\
If $N \geq d+2$ or if $N < d+2$ and $N$ does not
divide $d+1$, then the sphere $\sph^{d}$ cannot be embedded as a Lagrangian
submanifold of $W$, except possibly if $N=4$ and 
$d \equiv 2[4]$.
\end{thm}

\begin{proof}
As in the proof of Theorem~\ref{thm:nsc-spcut}, $N$ 
divides $2N_{e} = 2N_{W}$.
As $d \geq 2$, the sphere $L=\sph^{d}$ is simply connected and hence
exact with Maslov number $N_{L} = 2N_{W} = 2N_{e} \geq 4$. 
By Seidel's periodicity Theorem (Theorem~\ref{thm:S-per}), the Floer 
homology of $L$ is well defined, graded by $\zz/N$ and $2$-periodic.

We then use the following lemma which gives us the Floer homology of the
sphere for some values of $N_W$.

\begin{lem}
\label{lem:CFsph}
Let $W$ be a  symplectic monotone manifold of dimension $2d$ 
($d~\geq~2$), with first Chern number $N_{W} \geq 1$, admitting a 
sphere $\sph^{d}$ as Lagrangian submanifold.\\
Then the Floer homology of the sphere $\sph^{d}$ is well defined and 
if $2 N_{W}$ does not divide $d+1$, it is equal to its ordinary cohomology:
$$HF(\sph^{d},\sph^{d})= \bigoplus_{k \in \zz} H^{k}(\sph^{d};\zz/2)$$
\end{lem}

If $N$ does not divide $d+1$, then $2N_{W}$ cannot divide $d+1$, 
and by Lemma~\ref{lem:CFsph}, the Floer homology is isomorphic to the 
ordinary cohomology. This is in contradiction with the $2$-periodicity, 
except in the case $N=4$ and $d \equiv 2[4]$. 
\end{proof}

\textit{Proof of Lemma \ref{lem:CFsph}.}
As the sphere $L = \sph^{d}$ is simply connected for $d \geq 2$, it is
monotone and $N_L = 2 N_W \geq 2$. Thus, the Floer cohomology of $L$
is well defined.\\

If $N_L \geq d+2$, we can directly apply directly Oh's theorem
(Theorem \ref{thm:Oh}) and the Floer cohomology of the sphere is equal to
its ordinary cohomology.\\

To extend this result in the case $N_{L} < d+2$, we look precisely 
at the proof of the theorem of Oh. 

First, Oh defines a local Floer homology which is isomorphic to the 
Morse homology of $L$: given a Darboux neighbourhood $U$ of $L$ and an isotopy
$\Phi$ sufficiently close to the identity so that $\Phi_1(L)$ lies in
$U$, the local Floer homology is defined with the same complex than
the ``global'' Floer homology, but the differential counts only
the solutions of Floer's equation that lie in $U$. 

\begin{figure}[htbp]
  \begin{center}
   \psfrag{L}{$L$}
   \psfrag{Phi1(L)}{$\Phi_{1}(L)$}
   \psfrag{U}{$U$}
   \psfrag{x}{$x$}
   \psfrag{y}{$y$}
   \psfrag{u}{$u$}
   \includegraphics{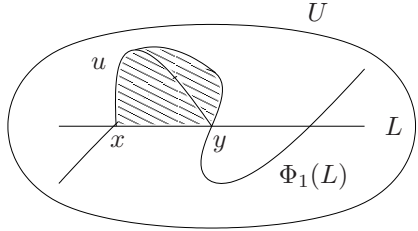}
    \caption{A solution $u$ that does not leave the neighbourhood $U$}
    \label{fig:reste}
  \end{center}
\end{figure}

We can restrict ourselves to a Hamiltonian isotopy constructed with the
help of a Morse function $f: L \rightarrow \rr$ on $L$: as $U$ is
symplectomorphic to a neighbourhood of the zero section in $T^{\ast}
L$, we define an Hamiltonian on $U$ by defining it on $T^{\ast} L$ by
the composition of $f$ and the projection $\pi$ of the cotangent
bundle: 
$$\begin{array}{cccc}
H: & T^{\ast} L & \rightarrow & \rr \\
   &   p        &   \mapsto   & f(\pi(p))
\end{array}$$
If $f$ is sufficiently $C^{1}$-small, the Hamiltonian isotopy defined
by $H$ is sufficiently close to the identity and if $f$ is
sufficiently $C^{2}$-small, we can even define an isomorphism between the
local Floer homology and the Morse cohomology (see \cite{MR1389956}
and \cite{MR1001276}).

We can then extend $H$ to the whole manifold $W$ by setting it equal
to zero out of (a neighbourhood of) $U$. Now, in order to prove the
same result for the (``global'') Floer homology (and thus the
theorem of Oh), it remains to see
that when $\Phi$ is close enough to the identity ($f$ is enough
$C^{2}$-small), all the solutions of Floer's equation that appear 
in the definition of the differential (i.e. the $J$-holomorphic strips 
between $L$ and $\Phi_{1}(L)$ with Maslov class equal to $1$) stay in $U$.\\

Assume that a solution $u$ ``gets out'' of $U$.

\begin{figure}[htbp]
  \begin{center}
   \psfrag{L}{$L$}
   \psfrag{Phi1(L)}{$\Phi_{1}(L)$}
   \psfrag{U}{$U$}
   \psfrag{u}{$u$}
   \psfrag{v}{$v$}
   \psfrag{x}{$x$}
   \psfrag{y}{$y$}
   \includegraphics{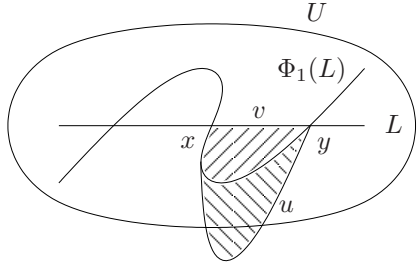}
    \caption{A solution $u$ that leaves $U$}
    \label{fig:sort}
  \end{center}
\end{figure}

We can glue to $u$ a small strip $v : [0,1] \times [1,2] \rightarrow U
\subset W$ in $U$ between $L$ and $\Phi_1(L)$ so that we get a disc
$w: [0,1] \times [0,2] \rightarrow W$ with boundary on $L$. 
The fact that the solution $u$ does not stay in $U$ implies (see
Proposition 4.1 of \cite{MR1389956}) that the symplectic area $\int
w^\star \omega \;$ of $w$ is positive. The monotonicity of $L$ then
implies that $\mu_L(w) > 0$.

But the Maslov classes of the strips satisfy:
$$\mu_L(w) = \mu_u(x,y) - \mu_v(x,y) = 1 - \mu_v(x,y)$$
and since $v$ is a small strip between $L$ and $\Phi_1(L)$, its Maslov
index is (see also \cite{MR926533}): 
$$\mu_v(x,y) = \mbox{ ind}_f(y) - \mbox{ ind}_f(x) \geq - \dim(L),$$
where $\mbox{ ind}_f(x)$ is the Morse index of $x$ for the function
$f$.

Hence, $\mu_L(w) = \mu_u(x,y) - \mu_v(x,y) = 1 - ( \mbox{ ind}_f(y) -
\mbox{ind}_f(x)) \leq 1 + d$.

But $\mu_L(w) > 0$ and by definition of the Maslov number, $N_L$
divides $\mu_L(w)$. We get here the (i) of Oh's Theorem: 
if $N_L \geq d+2$, all the solutions stay in $U$ and the Floer
homology is the Morse cohomology.\\

Here we are interested in the case when $L=\sph^{d}$ and we can choose
for the Morse function $f$ on $L$ the ``height'' function which has
only two critical points: a minimum of Morse index $0$ and a maximum
of Morse index $d$.

In this case, in order to have $\mu_L(w) >0$ and $N_L$ dividing
$\mu_L(w)$, it is necessary that $\mu_L(w) = 1+d$.

In consequence, if $N_{L} < d+2$ and $N_{L}$ does not
divide $d+1$, all the solutions stay in $U$, and the Floer homology 
is equal to the local homology, which is isomorphic to the ordinary
cohomology. \hfill \qedsymbol\\

\subsection{Monotone embeddings of tori}

The case of a monotone Lagrangian torus in a symplectic cut gives 
one more example of monotone Lagrangian torus with Maslov number 
equal to $2$. This is a result analogous 
to the result on monotone tori in $\cc^n$ (see \cite{LevMul, FO3}).

 \begin{thm} 
    \label{thm:tore-mondsspcut}
    Let $V^{d} \rightarrow B^{d-1}$ be a principal circle bundle with 
    basis a compact simply-connected manifold $B$ and non-vanishing Euler 
    class.\\
    Denote by $W$ a monotone symplectic cut of $\cotg$ (that is a cut at a 
    negative level).\\
    Assume that a torus $i: \TT^{d} \hookrightarrow W$ is embedded in $W$
    as a monotone Lagrangian submanifold with non trivial Maslov class.\\
    Then the Maslov number of $\TT$ is equal to $2$.
  \end{thm}

\begin{proof}
By \ref{prop:mon-W}, $W$ is monotone and its first Chern number is 
$N_{W}=N_{e}$ where $N_{e}$ is, as before, the Euler number of the 
circle bundle $V^{d} \rightarrow B^{d-1}$.

The torus being oriented, its Maslov number $N$ is even and greater or 
equal to $2$. Moreover, as $\TT$ is assumed to be monotone in $W$, the Floer 
homology of the torus $\TT$ is well defined. By the periodicity theorem of 
Seidel (Theorem~\ref{thm:S-per}), this homology $HF(\TT,\TT)$ is (absolutely) 
graded by $\zz/N$ and is $2$-periodic for this grading.

Let us assume that $N \geq 3$. We will get a contradiction by 
studying Oh's spectral sequence \cite{MR1389956}, and more 
precisely the description of this spectral sequence given by Biran.

Let $A = \zz/2[T,T^{-1}]$ be the algebra of Laurent polynomials over 
$\zz/2$ in the variable $T$. We define the degree of $T$ to be $N_L$.
Then $\displaystyle A = \bigoplus_{i \in \zz} A^i$ where  
$A^i = \zz/2 T^{i/N_L}$ if $i$ is a multiple of $N_L$ and 
$A^i = \{0\}$ otherwise.

\begin{thm}[Biran \cite{Birannew}] \label{thm:Biran}
Let $(M,\omega)$ be a geometrically bounded symplectic manifold and 
$L$ be a monotone Lagrangian submanifold of $(M,\omega)$ with 
Maslov number $N_{L} \geq 2$.\\
There exists a spectral sequence $\{E_{r}^{p,q},d_{r}\}$ satisfying 
the following properties:
\begin{description}
\item[(1)] For every $r \geq 0$, $E_{r}^{p,q}$ can be written  
  $E_{r}^{p,q}=V_{r}^{p,q} \otimes A^{pN_{L}}$ with  
  $d_{r}=\delta_{r} \otimes \tau_{r}$, where $\tau_{r}$ is the 
  multiplication by $T^{r}$, $V_{r}^{p,q}$ are vector spaces over 
  $\zz/2$, 
  $\delta_{r}: V_{r}^{p,q} \rightarrow V_{r}^{p+r,q-r+1}$ are  
  homomorphisms and satisfy $\delta_{r} \circ \delta_{r} = 0$. 
  Moreover,  
  $$V_{r+1}^{p,q}= 
\frac{\ker (\delta_{r}: V_{r}^{p,q} \rightarrow V_{r}^{p+r,q-r+1})}{\im (\delta_{r}: V_{r}^{p-r,q+r-1} \rightarrow V_{r}^{p,q} )};$$
\item[(2)] For $r=1$, $V_{1}^{p,q}=H^{p+q-pN_{L}}(L;\zz/2)$;
\item[(3)] $\{E_{r}^{p,q},d_{r}\}$ collapses at the $\nu + 1$ step, 
  where $\nu = [\frac{\dim(L)+1}{N_L}]$ and the spectral sequence converges to
  $HF(L,L)$, that is
  $$\bigoplus_{p+q=\ell} E_{\infty}^{p,q} \cong 
HF^{\ell(\modulo N_{L})}(L,L);$$
\item[(4)] For every $p \in \zz$, $$\bigoplus_{q \in \zz}
  E_{\infty}^{p,q} \cong HF(L,L).$$ 
\end{description}
\end{thm}

Moreover, Buhovski \cite{Lev} has endowed this spectral sequence of Biran of a 
multiplicative structure:

\begin{thm}[Buhovski, \cite{Lev}]
The spectral sequence of Theorem~\ref{thm:Biran} can be endowed with a multiplicative 
structure, that is, for every $r \geq 1$, there exists a product   
  $$m_{r}: E_{r}^{p,q} \otimes E_{r}^{p',q'} \longrightarrow
  E_{r}^{p+p',q+q'}$$
such that the differential $d_{r}$ satisfy the Leibniz rule with respect to 
this product. Moreover, this structure coincide on $E_{1}$ with the usual cup 
product on the cohomology $H^{*}(L;\zz/2)$.
\end{thm}

In the case of the torus $\TT$ in $W$, the differential at the first page is:
$$\delta_{1}: H^{p+q-pN}(\TT;\zz/2) \longrightarrow
H^{p+1+q-(p+1)N}(\TT;\zz/2).$$

In particular, if $p=0$ and $q=1$, 
$\delta_{1}: H^{1}(\TT;\zz/2) \rightarrow H^{2-N}(\TT;\zz/2)$ and hence,
if $N \geq 3$, $\delta_{1}$ vanishes on $V_{1}^{0,1}$. The cohomology of 
the torus is generated as a ring by the elements of 
$H^{1}(\TT;\zz/2)$, so that the Leibniz rule for $\delta_{1}$
with respect to the cup product implies that $\delta_{1}$ is zero. 

Consequently, 
$V_{2}^{p,q} = V_{1}^{p,q} = H^{p+q-pN_{\TT}}(\TT;\zz/2)$. 
But for $r \geq 2$, $\delta_{r}$ decreases the degree by $rN-1 > 3$. 
We get by induction $\delta_{r}=0$ and  
$E_{1}=E_{2}= ... = E_{r} = E_{\infty} = HF(\TT,\TT)$.

Eventually, $HF(\TT,\TT)=H^{*}(\TT;\zz/2)$. By Seidel's Theorem, we know that 
the Floer homology is graded 
by $\zz/N$ and $2$-periodic. This can be written:
for every $j \in \zz$, 
$$\bigoplus_{k \in \zz} H^{j+kN}(\TT;\zz/2)=\bigoplus_{k \in \zz}
H^{j+2+kN}(\TT;\zz/2).$$

At the dimension level, this equality can be written: for every $j \in \zz$, 
\begin{equation*}
\label{eq:1}
\sum_{k \in \zz} \binom {d} {j+ kN} = 
\sum_{k \in \zz} \binom {d} {j+2+kN}
\end{equation*}
(with the convention 
$\displaystyle \binom {d}{k} = 0$ 
if $k \not \in \{0, \ldots , d\}$).

This means, denoting for $j \in \{0, \ldots , N-1\}$, 
$$S_{j}=\sum_{k \in \zz} \binom{d}{j+kN}$$
that :
$$
\left \{
\begin{array}{l}
S_{0}= S_{2}= \ldots = S_{N-2}\\
S_{1}= S_{3}= \ldots = S_{N-1}.
\end{array}
\right.$$

As $N$ is even, $N=2m$, and we sum the $S_j$ for the even indexes $j$: 
$$S_{0} + S_{2} + \cdots + S_{2m-2} = \sum_{\ell=0, \: \ell \pair}^{d}  
\binom{d}{\ell} = 2^{d-1}$$
and for the odd indexes $j$:
$$S_{1} + S_{3} + \cdots  + S_{2m-1} = \sum_{\ell=0, \: \ell \impair}^{d}  
\binom{d}{\ell} = 2^{d-1},$$
and we get that $m S_{0} = m S_{1} = 2^{d-1}$. 
All the $S_j$ are thus equal and for every $j$ in $\zz$, $N S_j = 2^d$.

Note that, for a fixed $N$, if the property
\begin{equation}
\label{eq:2}
\mbox{for every } j \in \zz, \, N \sum_{k \in \zz} \binom{d}{j+kN} = 2^{d}
\end{equation}
is verified for $d=d_{0} \geq 2$, then by induction, using Pascal's rule, 
the property is verified for every $d \geq d_{0}$.

In order to get a contradiction, we prove that for every $d_{0} \geq 2$, 
there exists $d \geq d_{0}$ such that $N S_{0} \neq 2^{d}$.
For this purpose, we give an other expression of~$S_{0}$. 

Let
$\zeta=\exp(\frac{2 i \pi}{N})$ be a primitive $N$-root of the unity.
Expanding $(1+\zeta^{k})^{d}$ for $k \in \{0, \ldots , N-1\}$, and then 
summing up using the equalities:
$$1 + \xi + \xi^{2} + \cdots + \xi^{2m-1} = 0$$
where $\xi$ is a $N$-root of the unity, we get:
\begin{eqnarray*}
  N S_{0} & = & \sum_{k = 0}^{N-1} (1+\zeta^{k})^{d} \\
          & = & 2^{d} + \sum_{k = 1}^{N-1} (1+\zeta^{k})^{d}\\
          & = & 2^{d} + \sum_{k = 1}^{N-1} \cos \left(\frac{k d \pi}{N} \right)
\left( \cos \left( \frac{k \pi}{N} \right) \right)^{d}.
\end{eqnarray*}
For $d_{0} \geq 2$, choose $d \geq d_{0}$ a multiple of $2N$. 
Then $\displaystyle \cos \left( \frac{k d \pi}{N} \right) = 1$.
As for $k \in \{1, \ldots , N-1\}$, 
$0 < \dfrac{\pi}{N} \leq \dfrac{k \pi}{N} \leq
\dfrac{(N-1) \pi}{N} < \dfrac{\pi}{2}$ and 
$\cos \left( \dfrac{k \pi}{N} \right) > 0$, this implies
$$\sum_{k = 1}^{N-1}  \cos \left( \frac{k d \pi}{N} \right) 
\left( \cos \left( \frac{k \pi}{N} \right) \right)^{d} > 0$$
and 
$$N S_{0} > 2^{d}.$$
The Floer homology of the torus is thus not $2$-periodic and hence 
$N=2$. \end{proof}

\begin{rem}
We can also prove this result using instead of Biran's spectral sequence 
endowed with Buhovski's multiplicative structure the dichotomy of Biran 
and Cornea \cite{Bir-Cor}: the Floer homology of a monotone Lagrangian 
torus in a 
symplectic manifold is either trivial or isomorphic to its ordinary 
cohomology and in the case the Floer homology is trivial, the Maslov 
number is $2$. The proof above gives that the Floer homology cannot be 
isomorphic to the ordinary cohomology, so that the Floer homology of the 
torus must be trivial and its Maslov number is $2$.
\end{rem}

\subsection{Monotone embedding of the product of two spheres}
\label{sec:prodsphere}

Other examples of manifolds for which the cup product structure is 
interesting are the product of spheres $\sph^\ell \times \sph^m$:

\begin{thm} 
Let $V^{d} \rightarrow B^{d-1}$ be a principal circle bundle with 
basis a compact simply connected manifold $B$ and non vanishing Euler 
class.\\
Denote by $W$ a monotone symplectic cut of $\cotg$.\\
Assume that a product of two spheres $\sph^\ell \times \sph^m$ (for  
$1 \leq \ell \leq m$, $\ell+m=d$) can be embedded in $W$ as a monotone 
Lagrangian submanifold and that its Maslov class is not zero.\\
Then its Maslov number is less then or equal to $m+1$ except 
when $\ell=1$ and $m=2$ or when $\ell=4$ and $m=6$.
\end{thm}

\begin{proof}
By Kunneth Theorem, the ordinary cohomology (with $\zz/2$ coefficients) 
of $L=\sph^\ell \times \sph^m$ is, when $\ell < m$: 
$$H^k (L;\zz/2) = \left \{
\begin{array}{cl}
\zz/2 & \mbox{ if } k = 0,\ell,m,\ell+m \\
 \{0\}& \mbox{ otherwise,}
\end{array}
\right.$$
and when $\ell=m$:
$$H^k (L;\zz/2) = \left \{
\begin{array}{cl}
\zz/2              & \mbox{ if } k = 0,2m \\
\zz/2 \oplus \zz/2 & \mbox{ if } k = m \\
 \{0\}& \mbox{ otherwise.}
\end{array}
\right.$$

As in the case of the monotone torus, the Floer homology of $L$ is 
well defined, absolutely graded by $\zz/N$, if $N$ is the Maslov number 
of $L$, and is $2$-periodic for this grading.

At the first page of the spectral sequence, for $p=0$, 
the differential 
$$\delta_{1}: H^{q}(\TT;\zz/2) \rightarrow
H^{q+1-N}(\TT;\zz/2)$$
is trivial if $N \geq m+2$, so that $\delta_{1}$ is 
zero on $V_{1}^{0,q}$ for $q \leq m$. 
We investigate if  
$\delta_1:V_{1}^{0,q} \rightarrow V_1^{1,q}$ can be  
non zero for $q=\ell+m$. But $H^{\ell+m}(\TT;\zz/2)$ is generated 
by $H^{\ell}(\TT;\zz/2)$ and $H^{m}(\TT;\zz/2)$, and as for the torus, 
the Leibniz rule for $\delta_{1}$ with respect to the cup product 
implies that $\delta_{1}$ is zero. Analogously, $\delta_{1}$ is zero 
on every $V_{1}^{p,q}$.

Consequently, 
$V_{2}^{p,q} = V_{1}^{p,q} = H^{p+q-pN_{\TT}}(\TT;\zz/2)$. 
But for $r \geq 2$, $\delta_{r}$ decreases the degree by $rN-1 > q+1$. 
We deduce by induction that $\delta_{r}=0$ and  
$E_{1}=E_{2}= ... = E_{r} = E_{\infty} = HF(L,L)$.

Eventually, $HF(L,L)=H^{*}(L;\zz/2)$. As this Floer homology must 
be $2$-periodic, if $\ell < m$ and $N=m+2$, this is possible only 
when $m=2$ (and $\ell=1$) or when $m=6$ and $\ell=4$.

If $\ell<m$ and $N>m+2$ or if $\ell=m$, the Floer homology cannot be 
$2$-periodic.
\end{proof}

\setcounter{subsection}{0}

\setcounter{figure}{0}
\renewcommand{\thefigure}{\thesubsection.\arabic{figure}}

\section*{Appendices}
\addcontentsline{toc}{section}{Appendices}

\renewcommand{\thesubsection}{\Alph{subsection}}

\subsection{Complement to proofs}

\subsubsection{Complement to the proof of Lemma~\ref{lem:pi2mon}}
\label{sec:Append-Constr-map-Phi}

We give here the construction of the map $\Phi: \pi_{1}(\sph^{1})
\rightarrow \pi_{2}(W,V)$ fitting in the following commutative
diagram:
$$
\xymatrix{ \pi_{2}(V) \ar[r] \ar[d]^-{=} & \pi_{2}(B) \ar[r]
  \ar[d]^-{\simeq} &
  \pi_{1}(\sph^{1}) \ar[r] & \pi_{1}(V) \ar[r] \ar[d]^-{=} &  \pi_{1}(B) \ar[r] \ar[d]^-{\simeq}  &  0\\
  \pi_{2}(V) \ar[r] & \pi_{2}(W) \ar[r] & \pi_{2}(W,V) \ar[r] &
  \pi_{1}(V) \ar[r] & \pi_{1}(W) \ar[r] & 0 }
$$

Denote by $x_{0}$ the base point for $\pi_{1}(V)$, $b_{0}$ its
projection on $B$ and consider the fibre of $b_{0} \in B \subset Q$
for the complex line bundle $W \rightarrow B$. Adapting the definition 
of the diffeomorphism between $\cc$ and the symplectic cut of
$T^{\ast}\sph^{1}$ in Example 0 of Section \ref{sec:SC-cotg}, we
can describe a symplectic diffeomorphism between $\cc$ and the fibre
of $b_{0}$:
$$
\begin{array}{cccc}
i_{\cc} : & \cc & \longrightarrow & W \\
          &  z  & \longmapsto     & \left [ { x_{0}, \left({\xi +
                  \frac{1}{2} |z|^{2}} \right) \alpha_{x_{0}}, z } \right].
\end{array}
$$
The image of the disc of $\cc$ centered at the origin and of radius
$\sqrt{-2 \xi}$ is a disc in $W$ with boundary in $V$, since $V$ is
embedded in $\cotg$ as the zero section and in $W$ by:
$$
\begin{array}{ccc}
 V & \longrightarrow & W \\
 x & \longmapsto     & [x,0, \sqrt{-2 \xi}],
\end{array}
$$
and if $z = \sqrt{-2 \xi} u$ with $u \in \sph^{1}$, then 
$i_{\cc}(z) = [x_{0}, 0 , \sqrt{-2 \xi} u] = [\bar{u} \cdot x_{0}, 0 ,
\sqrt{-2 \xi}]$ lies in the image of $V$ in $W$.

Notice also that the center of the disc is mapped on the image of
$b_{0}$ in $W$: $i_{\cc}(0)= [x_{0},\xi \alpha_{x_{0}}, 0]$ (indeed,
$B$ is embedded in $Q$ by $[x] \mapsto [x,\xi \alpha_{x}]$ and $Q$ is 
embedded in $W$ by $[x,\varphi] \mapsto [x,\varphi,0]$).\\

In order to define $\Phi$, it is enough to give the image of the 
generator defined as the class of the loop  
$t \mapsto e^{-2 i \pi t}$ in $\pi_{1}(\sph^{1})$. We map this 
generator on the disc  
$$
\begin{array}{cccc}
\phi: &  D(0,\sqrt{-2 \xi}) & \longrightarrow & W \\
      &     z               & \longmapsto     & \left [ { x_{0},
          \left({\xi + \frac{1}{2} |z|^{2}} \right) \alpha_{x_{0}}, z } \right]
\end{array}
$$
of $W$ with boundary in $V$ which is the image by $i_{\cc}$ of the
disc of radius $\sqrt{-2 \xi}$ of $\cc$.
Let us check that this map $\Phi$ makes the diagram commute.\\

The two ``side'' squares
$$
\xymatrix{
  \pi_{1}(V) \ar[r] \ar[d]^-{=} & \pi_{1}(B)\ar[d]^-{\simeq}  & \mbox{\! and \!}  & \pi_{2}(V) \ar[r] \ar[d]^-{=} &  \pi_{2}(B) \ar[d]^-{\simeq}  \\
  \pi_{1}(V) \ar[r] & \pi_{1}(W) & & \pi_{2}(V) \ar[r] & \pi_{2}(W) }
$$
are naturally commutative as the homotopy equivalence between $B$ and
$W$ is given by the projection  $q_{\xi} \circ p_{\xi}$.\\

The square
$$
\xymatrix{
  \pi_{1}(\sph^{1}) \ar[r] \ar[d]^-{\Phi} & \pi_{1}(V)  \ar[d]^-{=} \\
  \pi_{2}(W,V) \ar[r] & \pi_{1}(V) }
$$
is commutative by construction of $\Phi$. Indeed, the image of the
generator $t \mapsto e^{-2 i \pi t}$ of $\pi_{1}(\sph^{1})$ by the
inclusion of $\sph^1$ in the fibre of $x_{0}$ is 
$u \mapsto \bar{u} \cdot x_{0}$ and the image of the class of $\phi$
by the boundary map is the class of 
$u \mapsto [x_{0}, 0, \sqrt{-2 \xi} \; u] = [\bar{u} \cdot x_{0}, 0 , \sqrt{-2 \xi}]$ 
which describes the loop $u \mapsto \bar{u} \cdot x_{0}$ of $V$ 
in $W$.\\

Finally, the square 
$$
\xymatrix{
  \pi_{2}(B) \ar[r] \ar[d]^-{\simeq} & \pi_{1}(\sph^{1}) \ar[d]^-{\Phi}\\
  \pi_{2}(W) \ar[r] & \pi_{2}(W,V) }
$$
is also commutative. The long exact sequence associated with the 
bundle $V \rightarrow B$ corresponds to the exact sequence of the pair 
$(V,\sph^1)$ where $\pi_{i}(V,\sph^{1})$ and $\pi_{i}(B)$ are identified by
the isomorphism $\pi_{\star}$ and the map 
$\pi_{i}(V,\sph^{1}) \simeq \pi_{i}(B) \rightarrow \pi_{i-1}(\sph^{1})$ 
is the boundary map.

If $w: \sph^{2} \rightarrow B$ is a sphere of $B$ with
$w(1)=b_{0}$, we can lift it to a disc $v$ of $V$ with boundary in the
fibre of $b_{0}$ (hence $\pi(v)=w$) such that $v(1)=x_{0}$. The class 
of $v$ in $\pi_{2}(V,\sph^{1})$ is then the image of the class of $w$ 
under the isomorphism 
$\pi_{2}(B) \stackrel{\simeq}{\rightarrow} \pi_{2}(V,\sph^{1})$.
Hence the image of the class of $w$ by the map 
$\pi_{2}(B) \rightarrow \pi_{1}(\sph^{1})$ is the class of the boundary 
$\partial v$ of $v$ and its image by the composition 
$\pi_{2}(B) \rightarrow \pi_{1}(\sph^{1}) \rightarrow \pi_{2}(W,V)$ is
the disc $\Phi(\partial v)$ constructed from the loop $\partial v$ of
$\sph^{1}$ by extenting it to a disc in the fibre of $b_{0}$.

On the other hand, the image of the class of a sphere $w$ of $B$ by
the map $\pi_{2}(B) \rightarrow \pi_{2}(W)$ is the class of $w$ after 
the embedding of $B$ in $W$ and the image of the class of a sphere $w$ of $W$
by $\pi_{2}(W) \rightarrow \pi_{2}(W,V)$ is simply the class of $w$
considered as a disc with boundary in $V$.

In order to verify that the diagram commutes, it is enough to check 
that in $\pi_{2}(W,V)$, the class of this sphere $w$ is the class of 
$\Phi(\partial v)$.

As the disk $v$ (the lift of the sphere $w$ of $B$) lies in $V$, its
class in $\pi_{2}(W,V)$ is trivial. As a consequence, the class in
$\pi_{2}(W,V)$ of the disc obtained by sum  of $v$ and of
$\Phi(\partial v)$ and denoted $v \# \Phi(\partial v)$ is
$$[v \# \Phi(\partial v)]= [\Phi(\partial v)].$$
But $v$ and $\Phi(\partial v)$ have the same boundary, hence 
$v \# \Phi(\partial v)$ is a sphere representing an element of 
$\pi_{2}(W)$. In consequence, we already know that the element
$[\Phi(\partial v)]$ of $\pi_{2}(W,V)$ is the image of an element of
$\pi_{2}(W)$. We still have to prove that $v \# \Phi(\partial v)$ is
homotopic to the sphere $w$ in $W$. But for this purpose, it is enough
to see that their image by the projection $q_{\xi} \circ p_{\xi}$ 
(which gives the homotopy equivalence between $W$ and $B$) are 
homotopic in $B$.
As $v$ takes value in $V$, $q_{\xi} \circ p_{\xi}(v)=\pi(v)=w$ and 
$q_{\xi} \circ p_{\xi}(\Phi(\partial v))=b_{0}$ in $B$ (recall that 
$\Phi$ is constructed such that the image of a loop is in the 
fibre of $b_{0}$). 
Eventually, $q_{\xi} \circ p_{\xi}(v \# \Phi(\partial v)) = w \# b_{0}$ 
which is homotopic to $w$, so that the last square is commutative. \\

\subsubsection{Proof of Proposition \ref{prop:poids}}
\label{sec:Append-Proof-Prop-poids}

The manifold $W$ being connected, it is enough to check that if $x$
and $y$ are two critical points such that there exists a trajectory of 
the gradient (of the Hamiltonian $H$ associated with the Hamiltonian
circle action) connecting these points, then 
$$w(x) \equiv w(y) [N_W] $$
where $w(x)$ denotes the sum of weights at the critical point $x$.

If $x$ and $y$ are two such critical points, we can construct a 
(gradient) sphere in $W$ having these points as poles. We will get the
result by studying the relationship between the first Chern class of this
sphere and the weights of the action at the poles.

Let $\varphi^t$ be a gradient trajectory connecting these
two points $x$ and $y$. It means that $\varphi^t$ satisfies 
$$ \frac{d \varphi^t}{dt} = - \grad _{\varphi^t} H, $$
$$\lim_{t \rightarrow - \infty} \varphi^t = x \mbox { and }  \lim_{t
  \rightarrow + \infty} \varphi^t = y.$$
In order to build the gradient sphere, we consider the trajectories 
$\varphi^t_u$ obtained from $\varphi^t$ by action of $\sph^{1}$, namely if 
$u \in \sph^{1}$, then 
$$\varphi^t_u : t \mapsto u \cdot \varphi^t$$
Let us prove that, for $u \in \sph^{1}$, $\varphi^t_u$ is also a trajectory of
the gradient.

We have (still denoting by $\rho$ the action)
\begin{eqnarray*}
\frac{d \varphi^t_u}{dt} & = & T_{\varphi^t} \rho(u) \, \frac{d  \varphi^t}{dt} \\
                         & = & - T_{\varphi^t} \rho(u) \, \grad_{\varphi^t} H.
\end{eqnarray*}
To conclude, it is enough to check that 
$T_{\varphi^t} \rho(u) \, \grad_{\varphi^t} H =  \grad_{u \cdot \varphi^t} H$ 
or more generally that, for all point $z$ of $W$, 
$T_{z} \rho(u) \, \grad_{z} H =  \grad_{u \cdot z} H$. This is true if
we choose a metric $g$ on $W$ which is invariant for the action of
$\sph^1$.

We have moreover 
$\displaystyle \lim_{t \rightarrow - \infty} (u \cdot \varphi^t) = u
\cdot ( \lim_{t \rightarrow - \infty} \varphi^t ) = u \cdot x = x$ 
(as $x$ is a critical point of $H$, it is a fixed point of the action)
and similarly, 
$\displaystyle \lim_{t \rightarrow + \infty} u \cdot \varphi^t = y$.
 
We can extend the ``cylinder without boundary''
$$
\begin{array}{ccc}
\sph^1 \times \rr & \rightarrow & W \\
       (u,t)      & \mapsto     & u \cdot \varphi^t
\end{array}
$$
to $\sph^1 \times [-\infty,+\infty] \rightarrow W$ by the limits $x$ and
$y$ respectively, and then by quotient of the boundaries of the
resulting  ``cylinder with boundary'', we get a gradient sphere 
$w : \sph^2 \rightarrow W$.

\begin{figure}[htbp]
  \begin{center}
   \psfrag{x}{$x$}
   \psfrag{y}{$y$}
   \psfrag{phit}{$\varphi^{t}$}
   \psfrag{uphit}{$u \cdot \varphi^{t}$}
   \includegraphics{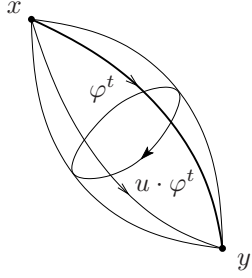}
    \caption{A gradient sphere between $x$ and $y$}
    \label{fig:spheregrad}
  \end{center}
\end{figure}

In order to get the first Chern class of the sphere $w$, we
consider the pullback by $w$ of the determinant bundle on $W$:
$$
\xymatrix{
w^\star \det(TW) \ar[r] \ar[d] & \det(TW) \ar[d] \\
   \sph^2 \ar[r]^-{w}             & W
         }
$$
The projection of the tangent bundle is equivariant under the action 
of $\sph^1$. It induces an equivariant projection 
$\det(W) \rightarrow W$.

Moreover, the map defined by:
$$
\begin{array}{ccc}
\sph^1 \times \rr & \rightarrow & W \\
    (u,t)      & \mapsto     & u \cdot \varphi^t
\end{array}
$$
is also equivariant if $\sph^1$ acts by multiplication on the first
factor of $\sph^1 \times \rr$. This action can be extended to  
$\sph^1 \times [-\infty,+\infty] \rightarrow W$ and then descends to the
quotient, so that the map $w$ and the projection of the line bundle 
$f^\star \det(TW) \rightarrow \sph^2$ are also equivariant.\\ 

But if we decompose the sphere $\sph^2$ in the union of the two open subset 
$\sph^2_+$ and $\sph^2_-$ obtained by removing from the sphere respectively the 
south pole and the north pole, we can trivialise this bundle on each of the two contractile 
open sets $\sph^2_+$ et $\sph^2_-$ : 
$$
\xymatrix{
f^\star \det(TW) \ar[rr]^-{\sim} \ar[dr] &         & \sph^2_+ \times \cc \simeq \cc \times \cc \ar[dl] \\
                                         &\sph^2_+ &
         }
$$
(and the same diagram holds for $\sph^2_-$).

If $w(y)$ (respectively $w(x)$) denotes the sum of the weights of 
the linearised action at the point $y$ (respectively at the point $x$), 
the action of $\sph^1$ on $\sph^2_+ \times \cc \simeq \cc \times \cc$ 
can be written: 
$$ u \cdot (z_1,\zeta_1) = (u z_1, u^{w(y)} \zeta_1)$$
and on $\sph^2_- \times \cc \simeq \cc \times \cc$ :
$$ u \cdot (z_2,\zeta_2) = (\bar{u} z_2, u^{w(x)} \zeta_2).$$
Indeed, the linearised action in the chart $\cc^n$ around $x$ is by
definition of the weights $(m_i(x))_{1 \leq i \leq n}$ at the point $x$: 
$$u \cdot \left( \begin{array}{c} v_1 \\ \vdots \\ v_n \end{array} \right) = 
\left( \begin{array}{ccc} u^{m_1(x)} &        & 0  \\ 
                                     & \ddots &    \\
                                0    &        & u^{m_n(x)} \end{array} \right)
\left( \begin{array}{c} v_1 \\ \vdots \\ v_n \end{array} \right)$$
and the determinant of the matrix is $u^{w(x)}$.\\
We can rebuild the complex line bundle $f^\star \det(TW) \rightarrow
\sph^2$ by gluing $\sph^2_+ \times \cc$ and $\sph^2_- \times \cc$ on their
common boundary $|z_1| = |z_2| = 1$. The gluing map can be expressed
in terms of the first Chern class:
$$\Phi (z_1,\zeta_1) = (\bar{z}_1, z_1^{c_1(f)} \zeta_1).$$
The action of $\sph^1$ must be compatible with the gluing, hence:
$$\Phi (u \cdot (z_1,\zeta_1)) = u \cdot \Phi (z_1,\zeta_1).$$
But 
\begin{eqnarray*}
\Phi (u \cdot (z_1,\zeta_1)) & = & \Phi (u z_1, u^{w(y)} \zeta_1) \\
                             & = & (\overline{u z_1}, (u z_1)^{c_1(f)} u^{w(y)} \zeta_1)
\end{eqnarray*}
and 
\begin{eqnarray*}
u \cdot \Phi (z_1,\zeta_1) & = & u \cdot (\bar{z}_1, z_1^{c_1(f)} \zeta_1) \\
                           & = & (\bar{u} \bar{z}_1, u^{w(x)} z_1^{c_1(f)} \zeta_1).
\end{eqnarray*}
Thus, we necessarily have: $c_1(f) + w(y) = w(x)$ and as $N_W$ divides 
$c_1(f)$, $w(x) \equiv w(y)\; [N_W]$. \hfil \qed

\subsection{Symplectic cut with semi-free circle action}
\label{sec:Append-Semilibre}

We assume in this section that the action of $\sph^1$ on $V$ is semi-free.
It induces as in the case of the free circle action on $V$ a Hamiltonian 
action on $\cotg$ of Hamiltonian:
$$H(x,\varphi)=\langle \varphi , X(x) \rangle.$$
The fixed points of the induced action on $\cotg$ sit in the $0$-level 
of $H$. All the other levels of $H$ are regular and the action of 
$\sph^1$ on these levels is free.
Consequently, the non-zero levels of $\tH$ are also regular and 
we can carry out the symplectic cutting at these levels.  
This time again, we will be interested in the symplectic cut 
at a negative level. 

In the following, we assume that $\xi$ is a negative real 
number. We want to know if in this case, the symplectic manifold 
(still denoted $W_\xi$ or $W$ for simplicity) is monotone.

\subsubsection{Sum of weights}
\label{sec:sompoids}

The fixed points of the circle action on the symplectic cut 
described at Section~\ref{sec:Seidel} are now not only the points 
of the symplectic submanifold $Q_\xi$ but also the fixed points 
in $H^{-1}(]\xi;+\infty[)$. Moreover, the fixed points for the action 
of the circle on the symplectic cut on the open 
set~$H^{-1}(]\xi;+\infty[)$ are the same as the fixed points for 
the action of the circle on $\cotg$. As we noticed above, 
these fixed points belong to $H^{-1}(0)$, and hence are embedded 
in the symplectic cut at the negative level $\xi$. 

In particular, any fixed point of $V$ is after the canonical embedding of $V$ 
in the symplectic cut a fixed point of the circle action 
on $W_\xi$. We prove that the sum of weights in these points is zero. 
It is enough to prove this fact for the circle action on the cotangent 
bundle.

\begin{lem}
The sum of the weights in a fixed point of the zero section $V$ for the action 
of the circle on the cotangent bundle is zero.
\end{lem}

\begin{proof}
We still denote by $\rho$ the action of $\sph^1$ on $V$ and we denote by 
$\trho$ the action induced on $\cotg$:
$$
\begin{array}{ccc}
\sph^1 \times (\cotg \times \cc) & \longrightarrow & \cotg \times \cc \\
       (u,(x,\varphi,z))         & \longmapsto     & (u \cdot x, u \cdot \varphi, \bar{u} z).
\end{array}
$$

In a fixed point $x$ of $V$, for every $u \in \sph^1$, 
$$T_x \rho(u): T_x V \rightarrow T_x V$$ 
is an endomorphism of $T_x V$. We choose again an invariant metric $g$ on $V$.
This means that $T_x \rho(u)$ is an isometry for the scalar product on $T_x V$ 
given by $g_x$. Denote $A(u)$ its matrix in an orthonormal basis for the 
scalar product on $T_x V$.

We know that the circle action on $\cotg$ is Hamiltonian, in particular 
symplectic. By definition, the sum of the weights in a fixed point 
is the sum of the exponents of the variable $u \in \sph^{1}$
when we write the matrix of the action in a common basis of 
diagonalisation.
In order to prove that the sum of weights is zero, it is enough to prove 
that the determinant of this matrix is equal to $1$ for every~$u$.

We express the circle action matrix on $\cotg$ in terms of the action 
matrix on~$V$. On the zero section of the cotangent bundle, the tangent 
space at $x$ can be decomposed in:
$$T_{(x,0)}(\cotg) \simeq T_x V \oplus T_x^\ast V.$$
In the basis obtained by considering the union of the orthonormal 
basis on $T_x V$ used to define $A(u)$ and its dual basis, the matrix 
of $T_x \trho(u)$ is:
$$\tA (u) =
\left( \begin{array}{cc} A(u) & 0 \\ 
                     0 & ^t A(u)^{-1} \end{array} \right)
= \left( \begin{array}{cc} A(u) & 0 \\ 
                              0 & A(u) \end{array} \right).
$$
Hence the matrix $\tA (u)$ is an element of the orthogonal matrices 
$O(2d)$ (as $A(u) \in O(d)$) and of the symplectic matrices $Sp(2d)$ 
(as the action is symplectic).

Moreover, we can choose for the almost complex structure 
$J_g$ on the cotangent bundle the one defined thanks to a 
Levi-Civit\`a connection for the metric on $V$ (as in 
\cite{MR1001276}, \cite{MR1389956}). This connection maps the 
vertical tangent vectors on the horizontal tangent vectors.
In particular, for the decomposition of the tangent space at $x$ 
above, if $v \in T_x V$, $J_g(v)=g(v,\cdot)$, so that the matrix 
$\tA(u)$ of $O(2d) \cap Sp(2d) \simeq U(d)$
correspond to the matrix $A(u)$ in $U(d)$.

We already know that $\det(A(u)) \in \{-1;1\}$ as $A(u)$ is orthogonal.
Moreover, the map $u \mapsto \det(A(u))$ is continuous and 
$\det(A(0))=\det(\Id)=1$, so that this map is a constant map 
equal to $1$ and the sum of the weights is equal to~$0$. \end{proof}

\subsubsection{Monotony of the symplectic cut}
\label{sec:topo}

In order to study the symplectic cut in the case of a semi-free 
circle action, we will use the manifold $\vreg$ obtained by 
taking away of $V$ the fixed points. The circle acts freely 
on $\vreg$ and we denote by $\breg = \vreg/{\sph^1}$ the quotient 
by the action.

As the action of the circle is semi-free, a consequence of the 
slice theorem (see for instance~\cite[Theorem I.2.1]{MR2091310}) 
is that the submanifold consisting of the fixed points is at 
least of codimension $2$ and $\vreg$ (and consequently also $\breg$) 
is path-connected.
 
The regular levels $H^{-1}(a)$ of $H$ (for $a \neq 0$) are, as in 
the case of the free circle action, isomorphic to vector bundles over 
$\vreg$ (they are isomorphic to the vector bundle 
$H^{-1}(0)_{|\vreg} \rightarrow \vreg$), hence are deformation 
retracts of $\vreg$.

Analogously, the quotients $Q_a = H^{-1}(a) / \sph^1$, 
for $a \neq 0$, are, as in the free action case, fibre bundle 
over $\breg$ isomorphic to $T^\ast \breg$ and hence have the 
homotopy type of $\breg$.\\

In order to prove that $W=W_\xi$ is monotone if $\xi < 0$, we study 
the group~$H_2(W)$. 

\begin{figure}[htbp]
  \begin{center}
   \psfrag{U}{$\cu$}
   \psfrag{V}{$\cv$}
   \psfrag{xi}{$\xi$}
   \psfrag{0}{$0$}
   \psfrag{H}{$H$}
   \psfrag{a}{$a$}
   \psfrag{b}{$b$}
   \psfrag{UintV}{$\cu \cap \cv$}
   \includegraphics{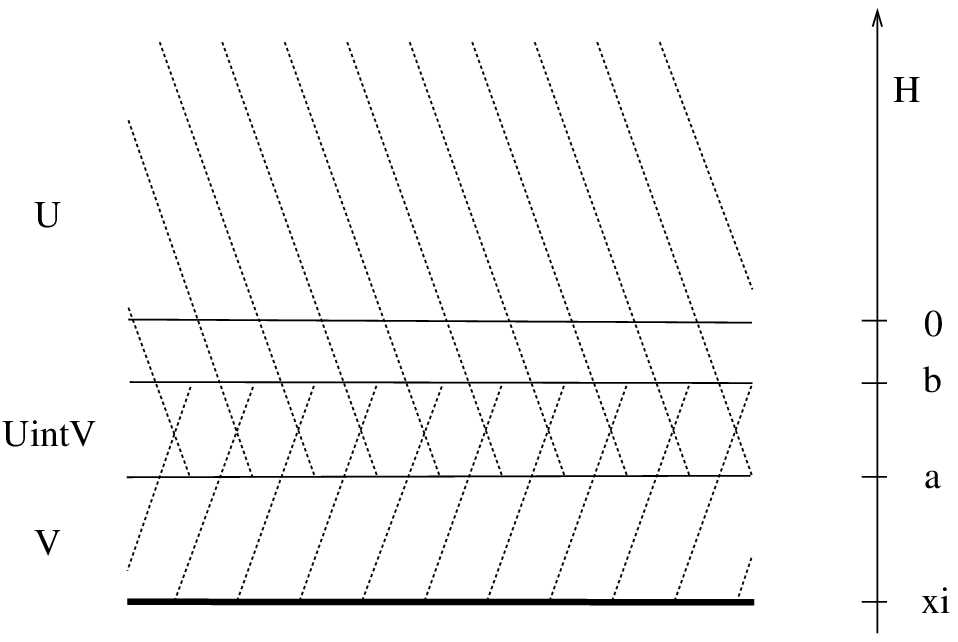}
   \caption{The union $W = \cu \cup \cv$}
   \label{fig:cotangentMV}
  \end{center}
\end{figure}

Let 
$$\cu = H^{-1}(a;+\infty) \mbox{ and } \cv =H^{-1}(\xi;b) \sqcup Q_\xi$$ 
for $\xi < a < b < 0$ (see Figure \ref{fig:cotangentMV}). 
These are two open sets of $W$ such that 
$\cu \cap \cv = H^{-1}(a;b)$, $V$ is a retract of $\cu$ (thanks to 
the restriction of the retraction of the cotangent bundle on its zero 
section), $Q_\xi$ is a retract of $\cv$ (it is a open disc fibre bundle 
over $Q_\xi$ with projection map the restriction of the one we had for 
$W \rightarrow Q$ in the case of the free action) and $\cu \cap \cv$ 
retracts on a regular level of $H$, thus has the homotopy type of 
$\vreg$.

We can determine $H_2(W)$ by using the Mayer-Vietoris exact sequence  
for the union $W=\cu \cup \cv$. 
The spectral sequence can be written here:
$$H_2(V) \oplus H_2(\breg) \stackrel{f_1}{\longrightarrow} H_2(W) 
\stackrel{f_2}{\longrightarrow} H_1(\vreg)
\stackrel{f_3}{\longrightarrow} H_1(V) \oplus H_1(\breg),$$
where $f_3$ is the map induced by the inclusion of $\vreg$ in $V$ 
and the projection of $\vreg$ on $\breg$.

If we have a section $s$ of the map 
$$\bar{f}_2: H_2(W) \longrightarrow \Ker{f_3}$$
induced by $f_2$, then by Mayer-Vietoris exact sequence,
 $H_2(W)$ is generated by the image in $H_2(W)$ of 
$H_2(V) \oplus H_2(\breg) \simeq H_2(V) \oplus H_2(Q)$ by the 
map $f_1$ and by the image 
of the section $s$.

We describe such a section $s$. We begin with the descripion of 
$\Ker{f_3}$. If $x \in V$ is a regular point in the neighbourhood 
of a codimension $2$ submanifold consisting of fixed points, 
then the orbit of this point $x$ under the circle action is a 
$1$-dimensional cycle in $\vreg$ whose image in $V$ by the inclusion 
and image in $\breg$ by projection are homologous to points.

We then prove that such cycles in the neighbourhood of the 
codimension $2$ fixed points submanifolds generate 
$\Ker{f_3}$. On the one hand, by the homotopy long exact sequence 
for the circle bundle $\vreg \rightarrow \breg$, the $1$-dimensional 
cycles of $\vreg$ which project on cycles homologous to a point 
in $\breg$ are homologous to orbits of points of $\vreg$ for the 
action of $\cerc$. On the other hand, if $\vvois$ denotes a 
neighbourhood given by the slice theorem of the fixed points 
submanifolds of~$V$, then the Mayer-Vietoris exact sequence 
for the union $V = \vreg \cup \vvois$ gives that the orbit of a 
point of $\vreg$ under the action of $\cerc$ is homologous 
to a point in $V$ if and only if this orbit is homologous to 
an orbit in $\vvois$.

A section $s$ can then be defined the following way. Let 
$x_0$ be a fixed point of the circle action on $V$. The point 
$[x_0,0,\sqrt{-2 \xi}] \in H^{-1}(0)$ is then a critical 
point of the Hamiltonian $h$ of the circle action on the 
symplectic cut (see Section~\ref{sec:Seidel}). 
Let $[x,\varphi,z]$ be a regular point in the unstable manifold of 
$[x_0,0,\sqrt{-2 \xi}]$ such that $x$ is in a neighbourhood 
(given by the slice theorem) of $x_0$. We assume that 
$a$ and $b$ are chosen such that $[x,\varphi,z]$ belongs to 
 $\cu \cap \cv$.

The trajectory $\gamma_t$ of the gradient of $h$ which goes through 
$[x,\varphi,z]$ at $t=0$ goes to $[x_0,0,\sqrt{-2 \xi}]$ as 
$t$ goes to $- \infty$. For $t > 0$, $\gamma_t$ lies in the 
disc bundle $\cv$ over $Q$ so that $\gamma_t$ goes, as $t$ goes 
to $+\infty$, to a fixed point $q$ which belongs to $Q$. Letting the 
circle act on this gradient trajectory, we get a gradient sphere 
with nord pole $x_0$, south pole $q$, whose ``equator'' 
$\cerc \cdot \gamma_0$ projects over $\cerc \cdot x$
in $\vreg$ which is a loop around the fixed point $x_0$.
This means that we can map any cycle around a codimension $2$ 
fixed points submanifold of $V$ to a gradient sphere in $W$ and 
this defines a section $s$ of $\bar{f}_2$.

Now, it is enough, in order to establish the monotonicity 
of $W$ to verify that the monotonicity relation is satisfy on the 
images of $H_2(\cu)$, $H_2(\cv)$ and on the gradient spheres 
in the image of the section $s$.

\begin{description}
\item[1)] The set $\cu$ is an open subset of $\cotg$ which is symplectically 
  embedded in~$W$. For any element $c$ of $H_2(\cu)$, 
  $c_1(W)(c)=0$ and $[\omega](c)=0$ (the form is exact on 
  this open set).

\item[2)] If we proceed to the symplectic cutting construction $W_{\reg}$
  at the level $\xi$ for the principal circle bundle 
  $\vreg \rightarrow \breg$, $\cv$ is a 
  disc bundle sitting inside the line bundle 
  $W_{\reg} \rightarrow Q_\xi$. We have proved that in the case of 
  a free action we have in $H^2(W_{\reg};\rr)$:
  $$[\omega](c) = - 2 \pi \xi c_1(W_{\reg})(c)= - 2 \pi \xi c_1(\cv)(c).$$

\item[3)] Let us look now at the monotonicity for the gradient spheres.
  Let $\gamma$ be a gradient trajectory of $h$ between a fixed point 
  $x_0$ of $V$ and a fixed point $q$ of $Q$. The circle action on 
  $\gamma$ defines a gradient sphere $\sph^1 \cdot \gamma$. Its area is 
  equal to:
  $$[\omega](\sph^1 \cdot \gamma) = 2 \pi \int_\gamma i_{X_h}
  \omega$$
  where $X_h$ is the Hamiltonian vector field associated to $h$ 
  (see for instance \cite[proof of Theorem VIII.1.1]{MR2091310}).
  Moreover, as $\gamma$ is a trajectory of the gradient of $h$,
  $$\frac{d \gamma}{dt}= - \grad \, h = J X_h.$$
  As a consequence, 
  $$\int_\gamma i_{X_h} \omega = \int \omega(X_h, J X_h) = 
  \int_\gamma \| X_h \|^2 = \int \left \| { \frac{d \gamma}{dt} } \right \| ^2$$
  is equal to the energy of this trajectory, that is:
  $$[\omega](\sph^1 \cdot \gamma) = 2 \pi (h(x_0)-h((x_{q}))= 2 \pi
  (- \xi).$$
  Besides, the Chern class of this sphere is the difference of the 
  sums of weights at $x_{0}$ and at $q$, that is
  $$c_1(\sph^1 \cdot \gamma) = w(q) - w(x_{0})= 1$$ 
  (see Sections~\ref{sec:Seidel}, \ref{sec:Append-Proof-Prop-poids} 
  and \ref{sec:sompoids}), so that the gradient 
  sphere satisfy  
  $$[\omega](c) = - 2 \pi \xi c_1(\cv)(c).$$
\end{description}

\bibliographystyle{plain}
\bibliography{Biblioanglaise}
\addcontentsline{toc}{section}{References}

\ \\
Agn\`es GADBLED\\
Mathematical Sciences Research Institute\\
17 Gauss Way\\
Berkeley, CA 94720-5070\\
e-mail: agnes.gadbled@gmail.com

\end{document}